\documentclass[11pt,a4paper]{article}
\usepackage{authblk}
\usepackage[utf8]{inputenc}
\usepackage[T1]{fontenc}
\usepackage{graphicx}
\usepackage{float}
\usepackage{subfigure}
\usepackage{multirow}
\usepackage{amssymb}
\usepackage{amsmath}
\usepackage{graphicx}

\usepackage[margin=1in]{geometry}
\usepackage{amsfonts}

\usepackage[numbers,sort&compress]{natbib}
\usepackage{hyperref}
\usepackage{caption}
\usepackage{subfigure}
\usepackage{mathrsfs}
\usepackage[title]{appendix}
\usepackage{xcolor}
\usepackage{textcomp}
\usepackage{manyfoot}
\usepackage{booktabs}
\usepackage{algorithm}
\usepackage{algorithmicx}
\usepackage{algpseudocode}
\usepackage{listings}
\allowdisplaybreaks
\newtheorem{thm}{Theorem}[section]

\newtheorem{deff}{Definition}[section]

\newtheorem{ex}{Example}[section]
\numberwithin{deff}{section}
\numberwithin{thm}{section}

\title{Stability and Bifurcation Analysis of Two-term Fractional Difference Equation}
\author{Janardhan Chevala}
\author{Sachin Bhalekar\footnote{Corresponding author email sachinbhalekar@uohyd.ac.in\\
 Email: janardhanch1992@gmail.com (J.C.)}}

\affil{School of Mathematics and Statistics, University of Hyderabad, Hyderabad, 500046 India.}

\date{}
\begin{document}
	
	\maketitle

	\begin{abstract}
	We consider the linear equation including two fractional order difference operators, viz. $\Delta^{\alpha}$ and $\Delta^{\beta}$, $0<\beta<\alpha \leq 1$. The sequence representation will be provided to find the solution in an easier way. The Z-transform will be used to find the boundary of the stable region in the complex plane. If the coefficient of the operator $\Delta^{\beta}$ is negative (near 0), then we observe that the boundary curve has multiple points generating multiple stability regions. We provide all possible bifurcations in terms of parameters. An ample number of examples will be provided to support the results.  
	\end{abstract}
	
%%%%%%%%%%%%%%%%%%%%%%%%%%%%%%%%%%%%%%%%%%%%%%%%%%%%%%%%%%%%%%%%%%%%%%%%%%%%%%%%%%%%%%%%%%%%%%%%%%
	
\section{Introduction} 
The mathematical modeling involves equations describing real-life phenomena \cite{aris1994mathematical,berry1995mathematical,gombert2000mathematical}. If the time variable involved in such equations is discrete, then the model is called a difference equation or a discrete dynamical system, or a map \cite{elaydi2005introduction}. A simple example of such an equation is the logistic equation $x(n+1)=r x(n) (1-x(n))$, $n=0,1,2,...$ describing the population of a species \cite{may1976simple}. Mathematical analysis of such models is an important aspect of difference equations. The equation $x(n+1)=f(x(n))$, $n=0,1,2,...$ where the state $x$ at the time $n+1$ depends on $x$ at previous time $n$ can be written as $\Delta x(n)=f(x(n))-x(n)$ by using a difference operator $\Delta$. We can also have more general cases like $x(n+1)=f(x(n), x(n-1), ..., x(n-k))$ which can be treated as delay difference equation \cite{joshi2024stability}. 
\par These equations can be used to design secure communication systems and in cryptography \cite{kotulski1999application,singh2010chaos,zhang2014novel}. Khrennikov \cite{khrennikov2004p} discussed applications of p-adic discrete dynamical systems to cognitive sciences and psychology. Huang and Wang  \cite{huang2015applications} established a model for population forecasting. Some more applications can be found in the following references list \cite{sandefur1990discrete,tu2012dynamical,bahi2013discrete,alsharawi2014theory}. 
\par Fractional calculus studies the operators of non-integer order \cite{podlubny1998fractional,kilbas2006theory}. The fractional differential equations are well studied in \cite{diethelm2002analysis,daftardar2004analysis,bhalekar2016stability,gupta2024fractional,bhalekar2025analysis}. These equations have applications in Bioengineering \cite{magin2012fractional}, Viscoelasticity \cite{mainardi2022fractional}, Physics \cite{hilfer2000applications} etc. We can also have a " fractional order difference operator", viz. $\Delta ^\alpha$, where $\alpha\in \mathbb{R}_{+}$ and $\Delta$ is forward difference operator to include the memory in the natural systems \cite{lubich1986discretized,miller1989fractional,atici2009initial,fulai2011existence,goodrich2015discrete,ferreira2022discrete}. The stability analysis of the fractional-order difference equation $\Delta^\alpha x(t)=f(x(t+\alpha-1))$, $0<\alpha<1$ is studied in \cite{abu2013asymptotic,vcermak2015explicit,bhalekar2022stabilitycomplexorder,bhalekar2022stabilitycoupledmap}. Bhalekar and Gade \cite{bhalekar2022stabilitycoupledmap} established stability conditions of synchronized fixed-point for linear fractional-order coupled map lattices. Stability of periodic map with period-2 in linear fractional difference equations is discussed by Bhalekar and Gade \cite{bhalekar2023fractional}. Joshi et al. proposed stability conditions of linear fractional order difference equations using feedback control \cite{joshi2023controlling} and with delay $\tau$ \cite{joshi2024stability}. In \cite{bhalekar2022stabilitycomplexorder}, Bhalekar et al. studied stability properties for complex order fractional difference equations.  In \cite{bhalekar2025dynamical}, a generalized logistic map with two parameters is proposed by Bhalekar et al.. They discussed stability bounds for each equilibrium point and also controlled chaos using delayed feedback.     
\par Bhalekar and Daftardar-Gejji proposed the theory of "multi-term" fractional differential equations in \cite{daftardar2008boundary,daftardar2008solving}. Cermak and Kisela \cite{vcermak2015stability} studied two-term fractional differential equations. Stability properties of a two-term fractional differential equation with delay are discussed by Bhalekar and Gupta \cite{bhalekar2024stability}. In \cite{vcermak2015asymptotic}, Cermak and Kisela discussed asymptotic stability conditions for a two-term linear fractional difference equation with backward difference operator $\nabla$. These equations are used in \cite{srivastava2010multi} to study a physiological system. Motivated by this work, we propose the two-term fractional order difference equation with Caputo-type operators. We consider two fractional order differences, viz. $\Delta^\alpha$ and $\Delta^\beta$ in the equation and study the stability properties. The rest of the article is organized as below:
\par Section 2 deals with the preliminaries. We describe the model in Section 3 and provide a sequence representation. Stability and bifurcation analysis are presented in Section 4. The conclusions are summarized in Section 5. 
%%%%%%%%%%%%%%%%%%%%%%%%%%%%%%%%%%%%%%%%%%%%%%%%%%%%%%%%%%%%%%%%%%%%%%%%%%%%%%%%%%%%%%%%%%%%%%%%%%
\section{Preliminaries} \label{prel}
In this section, we present some basic definitions and results.\\ 
Let $h>0$, \,$ a\in \mathbb{R}$,
$(h\mathbb{N})_a=\{a,a+h,a+2h,\dots\}$ and $\mathbb{N}_a=\{a,a+1,a+2,\dots\}$.
\begin{deff}(See\cite{bastos2011discrete,ferreira2011fractional,mozyrska2015transform})
	For a function $x:(h\mathbb{N})_a\longrightarrow \mathbb{R}$, the forward h-difference operator is defined as
	$$(\Delta_hx)(t)=\frac{x(t+h)-x(t)}{h},$$
	where $t\in(h\mathbb{N})_a$. \\
	Throughout this paper, we take $a=0$ and $h=1$. 
\end{deff}

\begin{deff}\cite{miller1989fractional,atici2007transform,atici2009initial,mozyrska2015transform}
	For a function $x:\mathbb{N}_0\longrightarrow \mathbb{R}$, the fractional sum of order $\beta >0$\, is given by 
	$$(\Delta^{-\beta} x)(t)=\frac{1}{\Gamma(\beta)} \sum_{s=0}^{t-\beta}{\frac{\Gamma(t-s)}{\Gamma(t-\beta-s+1)}}x(s),$$  
	where $t\in \mathbb{N}_{\beta}$.
\end{deff}

\begin{deff}\cite{fulai2011existence,mozyrska2015transform}  
	Let $\mu>0$ and $m-1<\mu<m$, where $m\in\mathbb{N}.$
	The $\mu$th fractional Caputo-like difference is defined as
	$$\Delta^\mu x(t)=\Delta^{-(m-\mu)}(\Delta^mx(t)),$$
	where $t\in\mathbb{N}_{m-\mu}$ and
	$$\Delta^mx(t)=\sum_{k=0}^m \left(\begin{array}{c}m\\k\end{array}\right)(-1)^{m-k}x(t+k).$$
\end{deff}

\begin{deff}\cite{mozyrska2015transform} 
	The Z-transform of a sequence $\{y(n)\}_{n=0}^\infty$ is a complex function given by $Y(z)=Z[y](z)=\sum_{k=0}^\infty y(k)z^{-k}$, where $z\in \mathbb{C}$ is a complex number for which the series converges absolutely.
\end{deff}
\begin{deff}\cite{mozyrska2015transform} 
	Let 
    \begin{eqnarray*}
		\Tilde{\phi}_\alpha(n) &=& \frac{\Gamma(n+\alpha)}{\Gamma(\alpha)\Gamma(n+1)} \nonumber \\
		&=& \left(\begin{array}{c} n+\alpha-1\\n \end{array} \right) = (-1)^n \left( \begin{array}{c}
			-\alpha \\ n \end{array} \right)
	\end{eqnarray*}
    be a family of binomial functions defined on $\mathbb{Z}$, parametrized by $\alpha$.\\
Then, 
	\begin{eqnarray*}
		Z(\Tilde{\phi}_\alpha(t))=\frac{1}{(1-z^{-1})^\alpha},\; |z|>1. \nonumber
	\end{eqnarray*}
\end{deff}
\begin{deff}\cite{mozyrska2015transform}
	The convolution $\phi\ast x$ of the functions $\phi$ and $x$ defined on $\mathbb{N}_0$ is defined as 
	\begin{eqnarray*}(\phi\ast x)(n)=\sum_{s=0}^n \phi(n-s)x(s)=\sum_{s=0}^n \phi(s) x(n-s).
	\end{eqnarray*} 
	Then, the Z-transform of this convolution is
	\begin{eqnarray*}
		Z(\phi\ast x)(n)=(Z(\phi)(n)) Z((x)(n)).
	\end{eqnarray*} 
\end{deff}

\begin{deff}\cite{elaydi2005systems,hirsch2013differential}
	A steady state solution or an equilibrium $x_*$ of the fractional difference equation $$\Delta^{\alpha} x(t)=f(x(t+\alpha-1))-x(t+\alpha-1),$$ where $0<\alpha \leq 1$, $t \in \mathbb{N}_{1-\alpha}$ is a real number satisfying $f(x_*) = x_*$.
\end{deff}

\begin{deff} \cite{elaydi2005systems,hirsch2013differential}
	An equilibrium $x_*$ is stable if for each  $\epsilon>0$, there exists $\delta>0$ such that $|x_0 -x_* | < \delta $ implies
	$|x(t) - x_* | < \epsilon$, $t=1,2,3,...$\\
	If $x_*$ is not stable, then it is unstable.
\end{deff}

\begin{deff} \cite{elaydi2005systems,hirsch2013differential}
	An equilibrium point $x_*$  is asymptotically stable  if it is stable and there exists $\delta>0$ such that $|x_0 -x_* | < \delta $ implies $ \lim_{t\to\infty}x(t) =x_*$.
\end{deff}
%%%%%%%%%%%%%%%%%%%%%%%%%%%%%%%%%%%%%%%%%%%%%%%%%%%%%%%%%%%%%%%%%%%%%%%%%%%%%%%%%%%%%%%%%%%%%%%%%%
\section{The model}

In this work, we consider the two-term linear fractional difference system  as
\begin{equation*} 
\Delta^{\alpha} x(t) + a \, \Delta^{\beta} x(t+\alpha-\beta) =f(x(t+\alpha-1))-x(t+\alpha-1), 
\end{equation*}
\begin{equation}
   x(0)=x_0, \label{1}
\end{equation}
 with $f(x)=bx$, where $0<\beta<\alpha  \leq 1$, $t \in \mathbb{N}_{1-\alpha}$ and $a\in \mathbb{R}$, $b \in \mathbb{C}$ are constants.\\
\par In this paper, we study the stability analysis of a two-term fractional difference equation (\ref{1}).\\

Now, we provide an equivalence between the initial value problem (\ref{1}) and a sequence $x(n)$. The sequence representation will be extremely useful while solving examples numerically.

\begin{thm}
The initial value problem (\ref{1}) is equivalent to the sequence representation 
\begin{equation}
\begin{split}
x(n)&=\left(\frac{\alpha+a \beta+(b-1)}{a+1} \right) x(n-1)\\
    & +\frac{1}{a+1} \sum_{s=0}^{n-2}\frac{1}{(n-s)!} \left( \frac{\alpha \Gamma(n-s-\alpha)}{\Gamma(1-\alpha)}+ \frac{a \beta\Gamma(n-s-\beta)}{\Gamma(1-\beta)}\right) x(s)
\end{split} \label{seqrep}
\end{equation}
where $n=1,2,3,...$.
\end{thm}
\textbf{Proof.} \quad
    Consider $ \Delta^{\alpha} x(t) + a \,  \Delta^{\beta} x(t+\alpha-\beta) =(b-1) \, x(t+\alpha-1)$. \\
By the def. of Caputo fractional difference, the above  equation is equivalent to\\

$\Delta^{-(1-\alpha)}(\Delta x(t))+a \, \Delta^{-(1-\beta)}(\Delta x(t+\alpha-\beta))=(b-1) \, x(t+\alpha-1) $.\\
Applying the def. of fractional sum,\\
\begin{equation}
\begin{array}{l}
\displaystyle
\frac{1}{\Gamma(1 - \alpha)} \left( 
    \sum_{s=0}^{t+\alpha} \frac{\Gamma(t+1-s)}{\Gamma(t+1-s+\alpha)} x(s)
    - \sum_{s=0}^{t-1+\alpha} \frac{\Gamma(t-s)}{\Gamma(t-s+\alpha)} x(s)
\right) \\[12pt]
\displaystyle
+ \frac{a}{\Gamma(1 - \beta)} \left( 
    \sum_{s=0}^{t+\alpha} \frac{\Gamma(t+1-s+\alpha-\beta)}{\Gamma(t+1-s+\alpha)} x(s)
    - \sum_{s=0}^{t-1+\alpha} \frac{\Gamma(t-s+\alpha-\beta)}{\Gamma(t-s+\alpha)} x(s)
\right)\\[12pt]
\displaystyle
= (b - 1)\, x(t + \alpha - 1). \label{step1}
\end{array}
\end{equation}
Putting $t=n-\alpha$ in the eq.(\ref{step1}), we get
\begin{equation}
\begin{array}{l}
\displaystyle
\frac{1}{\Gamma(1 - \alpha)} \left( 
    \sum_{s=0}^{n} \frac{\Gamma(n-s+1-\alpha)}{\Gamma(n-s+1)} x(s)
    - \sum_{s=0}^{n-1} \frac{\Gamma(n-s-\alpha)}{\Gamma(n-s)} x(s)
\right) \\[10pt]
\displaystyle
+ \frac{a}{\Gamma(1 - \beta)} \left( 
    \sum_{s=0}^{n} \frac{\Gamma(n-s+1-\beta)}{\Gamma(n-s+1)} x(s)
    - \sum_{s=0}^{n-1} \frac{\Gamma(n-s-\beta)}{\Gamma(n-s)} x(s)
\right) 
- (b - 1)\, x(n - 1) = 0.
\end{array}
\label{step2}
\end{equation}
Now collecting the coefficients of $x(s)$ (for s=0,1,...n) and simplifying eq.(\ref{step2}), we get
\begin{equation}
\begin{array}{l}
\displaystyle
(a + 1)\, x(n) - (\alpha + a\beta + b - 1)\, x(n - 1) \\[10pt]
\displaystyle
\quad -  
\sum_{s=0}^{n-2} \frac{1}{\Gamma(n - s + 1)} \left( 
\frac{\alpha\, \Gamma(n - s - \alpha)}{\Gamma(1 - \alpha)} 
- \frac{a\beta\, \Gamma(n - s - \beta)}{\Gamma(1 - \beta)} 
\right) x(s) = 0.
\end{array}
\label{step3}
\end{equation}

From this eq. (\ref{step3}), We get the required eq. (\ref{seqrep}).\\ 
The converse part can be proved in a similar way.\\
This proves the result.

%%%%%%%%%%%%%%%%%%%%%%%%%%%%%%%%%%%%%%%%%%%%%%%%%%%%%%%%%%%%%%%%%%%%%%%%%%%%%%%%%%%%%%%%%%%%%%%%%%
 \section{Stability and bifurcation analysis}
 \subsection*{Characteristic equation}
 Taking Z-transform on both sides of (\ref{1}), we get
 \begin{eqnarray*}
   p [(z-1)X(z)-zx(0)]+a \,  q [(z-1)X(z)-zx(0)]=(b-1)X(z),
 \end{eqnarray*}
 where $p=\left(\frac{z}{z-1}\right)^{1-\alpha}$, $q=\left(\frac{z}{z-1}\right)^{1-\beta}$ and X is the Z-transform of x.
 \begin{equation}
 X(z)\left[(z-1)p+  a (z-1) q-(b-1)\right]=zx(0)\left[p+a q\right].  \label{ztran}
 \end{equation}
 The characteristic equation of the system (\ref{1}) can be obtained by equating the coefficients of the term $X(z)$ in eq. (\ref{ztran}) to zero as  
 \begin{equation}
     z(1-z^{-1})^\alpha+a \, z(1-z^{-1})^\beta+1=b,\label{chareq}
 \end{equation}
 where the condition $|z|<1$ should be satisfied for the stability. Thus, the boundary of the stable region is given by $|z|=1$.\\
 Therefore, we can find the boundary of the stable region by setting $z=e^{i\theta}$ in equation (\ref{chareq}).\\ We get
 \begin{equation}
     e^{i\theta}(1-e^{-i\theta})^\alpha+a \, e^{i\theta}(1-e^{-i\theta})^\beta+1=b,\label{bval}
 \end{equation}
 where $\theta \in [0, 2 \pi]$.\\
 Simplify Eq. (\ref{bval}), we get
 \begin{eqnarray}
   b &=& 2^\alpha \left(\sin\frac{\theta}{2}\right)^\alpha e^{i [\frac{\alpha \pi}{2}+\theta(1-\frac{\alpha}{2})]}+a \,  2^\beta \left(\sin\frac{\theta}{2}\right)^\beta e^{i [\frac{\beta \pi}{2}+\theta(1-\frac{\beta}{2})]} +1.\label{simbval}
\end{eqnarray}
 \subsection*{Boundary curve}
 The parametric representation of the boundary curve after separating real and imaginary parts of Eq. (\ref{simbval}) is
 \begin{equation}
     \gamma(\theta)=\left( \gamma_1(\theta) , \gamma_2(\theta)\right),\label{boundarycurve}
 \end{equation}
 where 
\begin{eqnarray*}
    \gamma_1(\theta)=2^\alpha \left(\sin\frac{\theta}{2}\right)^\alpha \left(\cos\left(\frac{\alpha \pi}{2}+\theta(1-\frac{\alpha}{2}\right)\right)+a \, 2^\beta \left(\sin\frac{\theta}{2}\right)^\beta \left(\cos\left(\frac{\beta \pi}{2}+\theta(1-\frac{\beta}{2}\right)\right) +1 ,
    \end{eqnarray*}
  and $ \gamma_2(\theta)= 2^\alpha \left(\sin\frac{\theta}{2}\right)^\alpha \left(\sin\left(\frac{\alpha \pi}{2}+\theta(1-\frac{\alpha}{2}\right)\right)+a \, 2^\beta \left(\sin\frac{\theta}{2}\right)^\beta \left(\sin\left(\frac{\beta \pi}{2}+\theta(1-\frac{\beta}{2}\right)\right)$. \\
  Here $\gamma_1(\theta)$ and $\gamma_2(\theta)$  are real and imaginary parts of the expression $b$ in (\ref{simbval}).
 \par If the complex number $b$ lies inside the positively oriented simple closed curve $\gamma(\theta)$, then the system (\ref{1}) will be asymptotically stable.\\

\begin{figure}[H]
    \centering
    \subfigure[$\alpha=0.9, \, \beta=0.2$]{
\includegraphics[height=1.6in,width=1.6in,keepaspectratio]{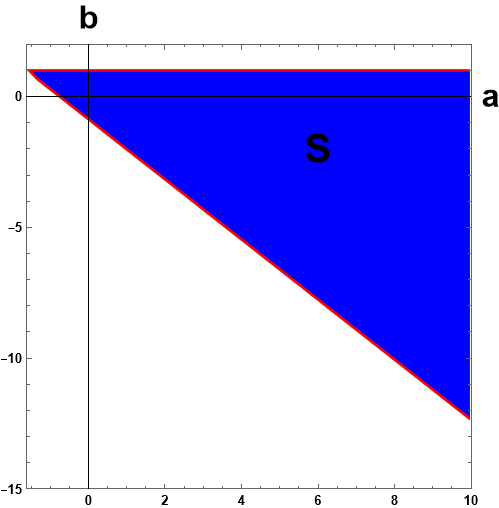}
        \label{reg1}
    } \hspace{0.1cm}
     \subfigure[$\alpha=0.9, \, 0.1 \leq \beta \leq 0.8$]{
\includegraphics[height=2.8in,width=2.8in,keepaspectratio]{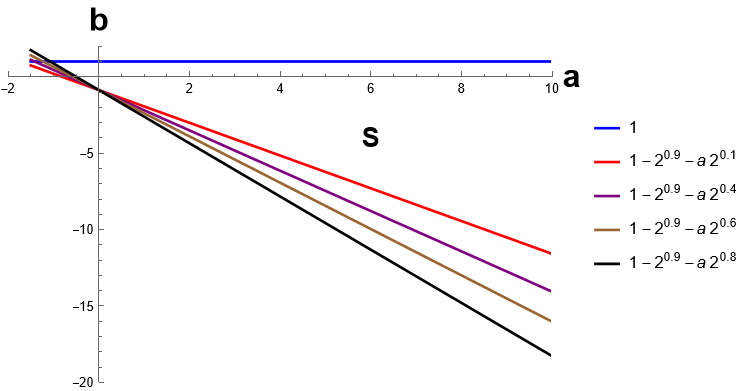}
        \label{reg2}
    } \hspace{0.1cm}
     \subfigure[$\alpha=0.7, \, 0.1 \leq \beta \leq 0.6$]{
\includegraphics[height=2.8in,width=2.8in,keepaspectratio]{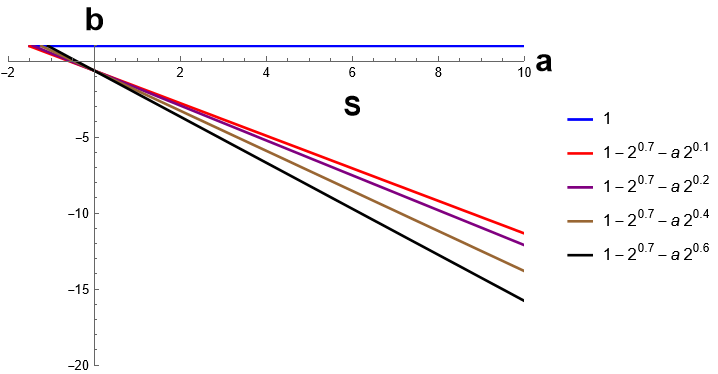}
        \label{reg3}
    } \hspace{0.1cm}
    \subfigure[$\beta=0.1, \, 0.2 \leq \alpha \leq 0.8$]{
\includegraphics[height=2.8in,width=2.8in,keepaspectratio]{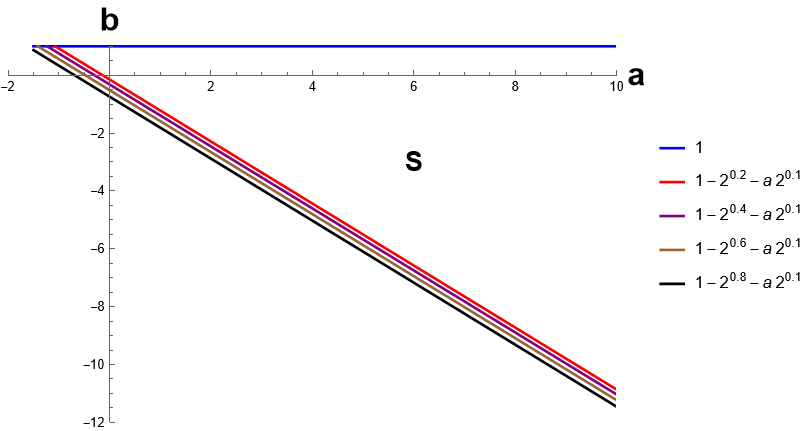}
        \label{reg4}
    }
    \caption{ Stable regions of system (\ref{1}) for different values of $\alpha$ and $\beta$ in $ab$-plane. In figures, S: stable.}
    \label{reg}
\end{figure}

\textbf{Note}: If $a$ and $b$ are reals, then the boundary of the stable region of the system (\ref{1}) can be obtained by substituting $z=\pm 1$ in the characteristic equation (\ref{chareq}). We get, the stable region of the system (\ref{1}) as $1-2^{\alpha}-a \, 2^{\beta}<b<1$. This region is bounded by the straight lines $b=1$ and $b=1-2^\alpha-a \, 2^\beta$ in the $ab-$plane. These lines intersect at $(-2^{\alpha-\beta},1)$ in the $ab-$plane. We sketch the stable regions of system (\ref{1}) in the $ab-$plane for various values of $\alpha$ and $\beta$ (
cf. Fig.\ref{reg} ). \\

\begin{figure}[H]
    \centering
    \subfigure[$b=-1.5$]{
\includegraphics[height=2.1in,width=2.1in,keepaspectratio]{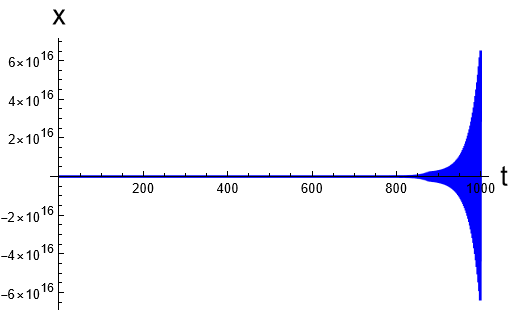}
        \label{re1a}
    } \hspace{0.3cm}
    \subfigure[$b=-1.3$]{
\includegraphics[height=2.1in,width=2.1in,keepaspectratio]{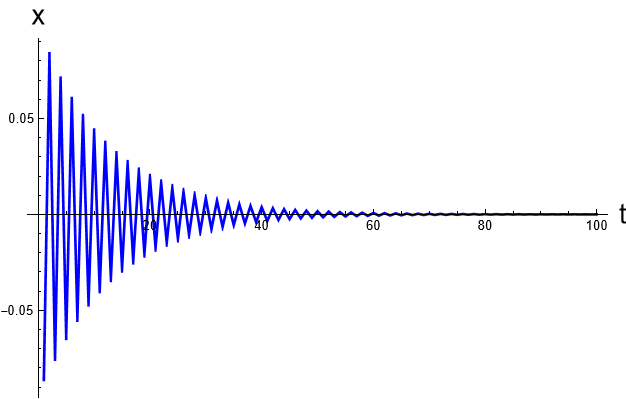}
        \label{re1b}
    } \hspace{0.3cm}
    \subfigure[$b=0.8$]{
\includegraphics[height=2.1in,width=2.1in,keepaspectratio]{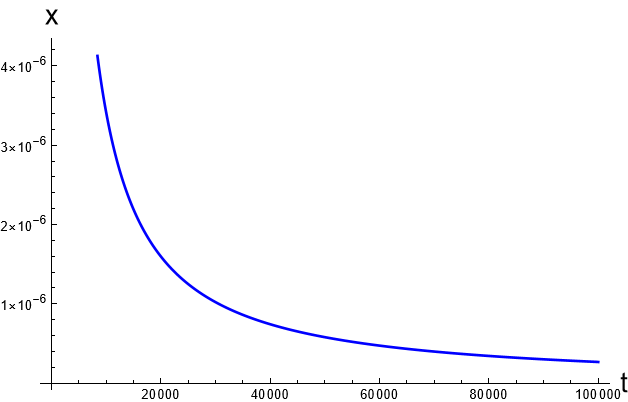}
        \label{re1c}
    } \hspace{0.3cm}
    \subfigure[$b=2.5$]{
\includegraphics[height=2.1in,width=2.1in,keepaspectratio]{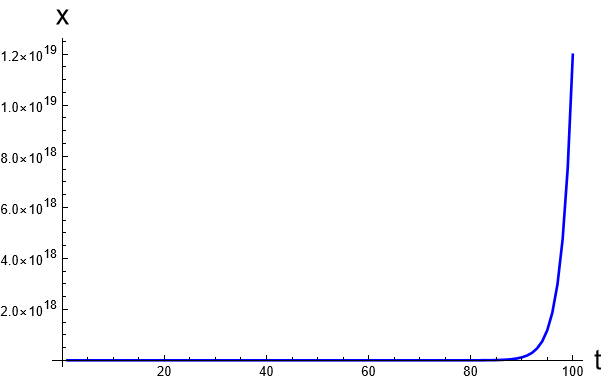}
        \label{re1d}
    }
    \caption{ The solutions of system (\ref{1}) for different values of $b$ within and outside the stable region.}
    \label{re}
\end{figure}

\begin{ex}
    Consider $\alpha=0.8$, $\beta=0.2$ and $a=0.6$. For $b \in (-1.43032,1)$, the system (\ref{1}) is stable in this case. Let $b=-1.3$ and $b=0.8$ be sample values in the stable region and the initial condition x(0)=0.1. Thus, the system (\ref{1}) converges to $0$ ( cf. Fig.\ref{re1b} and Fig.\ref{re1c} ). However, for $b=-1.5$ and $b=2.5$ outside the stable region on the left and right, respectively, the system (\ref{1}) is unbounded ( cf. Fig.\ref{re1a} and Fig.\ref{re1d} ).
\end{ex}

\par For $a=0$ in system (\ref{1}), we get classical case $\Delta^{\alpha} x(t) =f(x(t+\alpha-1))-x(t+\alpha-1)$, where the stable region in the complex plane is given by the cardioid as discussed in the literature ( see Fig.1 in \cite{bhalekar2022stabilitycoupledmap} ).\\
 
 \par
\begin{figure}[H]
    \centering
    \subfigure[$\alpha=0.4$(red), $0.6$(purple) and $0.8$(blue) resp.]{
\includegraphics[height=2.6in,width=2.6in,keepaspectratio]{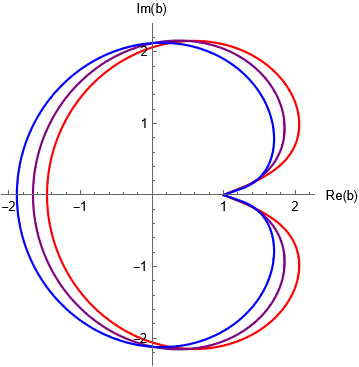}
        \label{beta_0.2}
    } \hspace{0.1cm}
    \subfigure[$\alpha=0.7$(red), $0.8$(black) and $0.9$(blue) resp.]{
\includegraphics[height=2.6in,width=2.6in,keepaspectratio]{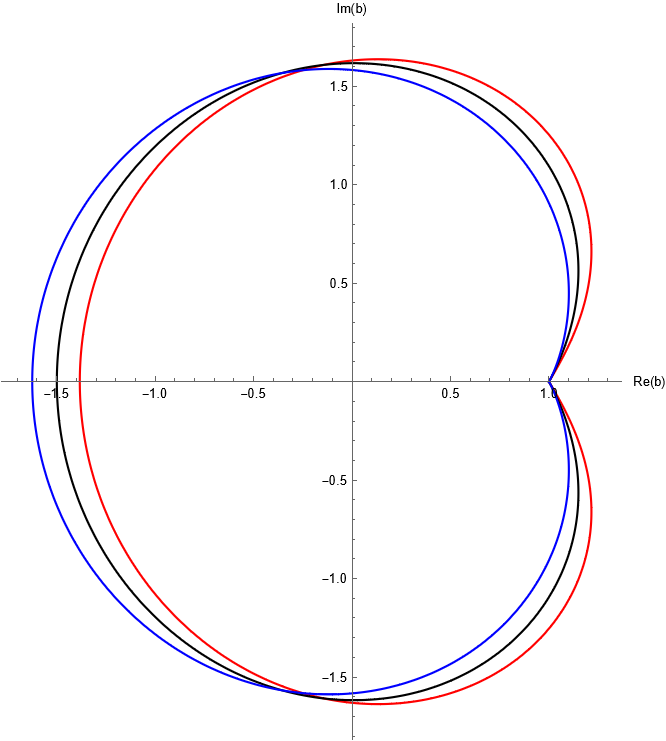}
        \label{beta_0.6}
    } \hspace{0.1cm}
    \subfigure[$\beta=0.2$(red), $0.4$(purple) and $0.6$(blue) resp.]{
\includegraphics[height=2.6in,width=2.6in,keepaspectratio]{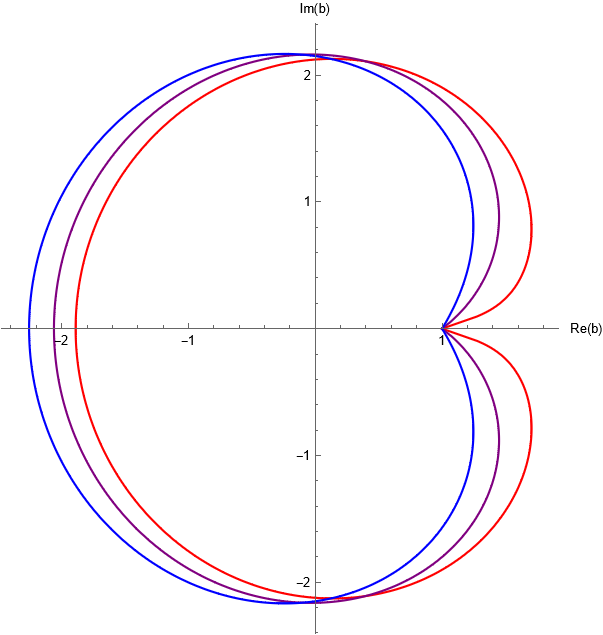}
        \label{alpha_0.8}
    } \hspace{0.1cm}
    \subfigure[$\beta=0.05$(red), $0.3$(black) and $0.45$(blue) resp.]{
\includegraphics[height=2.6in,width=2.6in,keepaspectratio]{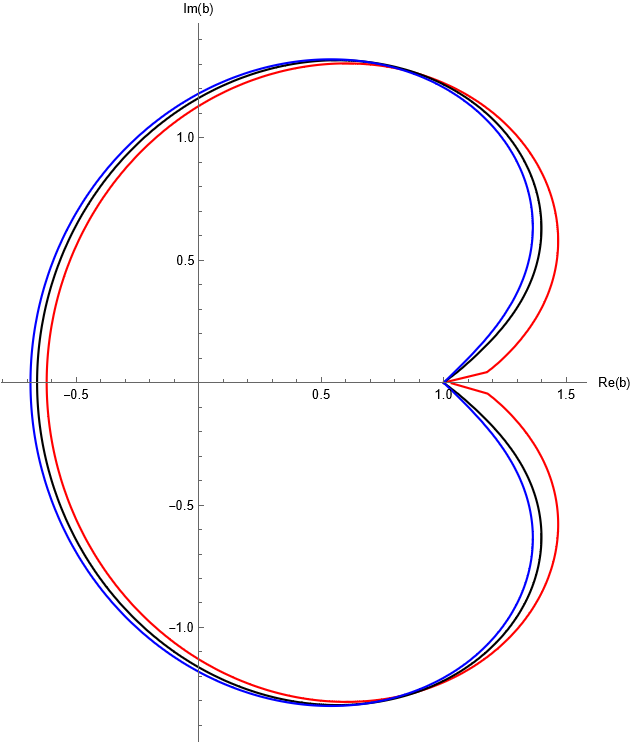}
        \label{alpha_0.5}
    }
    \caption{ Stable regions of system (\ref{1})  for fixed values $\beta=0.2$, $a=1$ in (a), $\beta=0.6$, $a=0.5$ in (b), $\alpha=0.8$, $a=1$ in (c), and $\alpha=0.5$, $a=0.2$ in (d) in a complex plane, respectively.}
    \label{alphabeta}
\end{figure}

 For $a>0$, the curve $\gamma(\theta)$ traces a cardioid with cusp at $b=1$. For fixed $\alpha, \beta$ as $a$ increases, the cusp of the cardioid at $b=1$ becomes smoother and vanishes eventually, resulting in a circular stable region. Here, the stable region increases on the left side, and it is symmetric about the x-axis. We sketch the stable regions of system (\ref{1}) for fixed $a$ and $\alpha$, and various values of $\beta$ ( see Figs. \ref{alpha_0.8} and \ref{alpha_0.5} ) and also sketch the stable regions of system (\ref{1}) for fixed $a$ and $\beta$, and various values of $\alpha$ ( see Figs. \ref{beta_0.2} and \ref{beta_0.6} ) in the complex plane.\\
               
\begin{figure}[H]
    \centering
    \subfigure[$b=0.5+i$]{
\includegraphics[height=2.5in,width=2.5in,keepaspectratio]{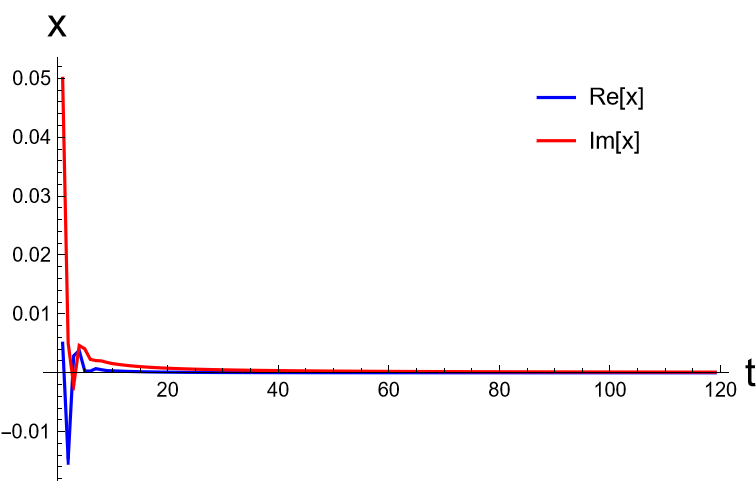}
        \label{pa1a}
    } \hspace{0.3cm}
    \subfigure[$b=-1.5+i$]{
\includegraphics[height=2.5in,width=2.5in,keepaspectratio]{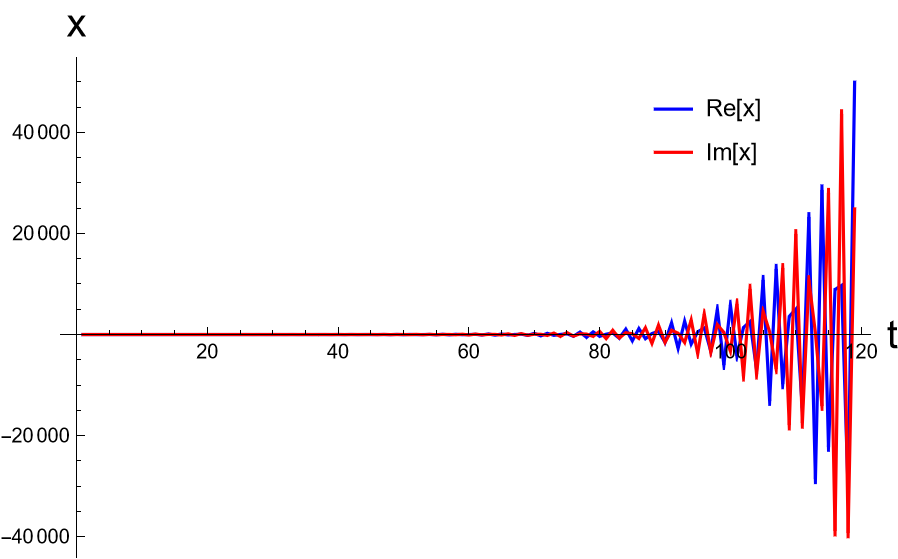}
        \label{pa1b}
    } \hspace{0.3cm}
    \subfigure[$b=2.5+i$]{
\includegraphics[height=2.5in,width=2.5in,keepaspectratio]{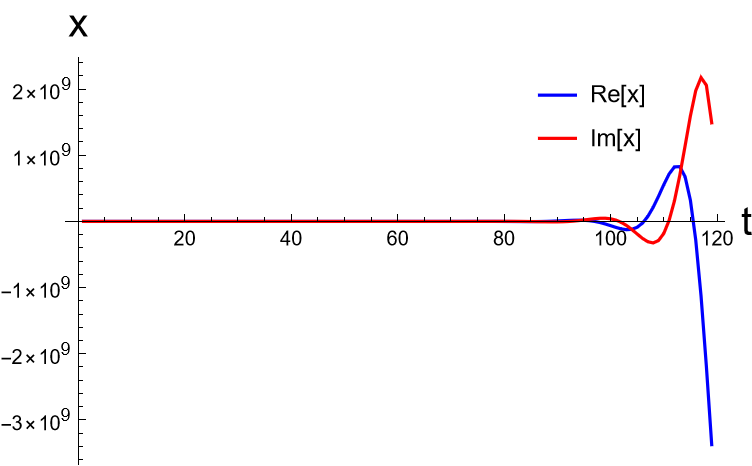}
        \label{pa1c}
    }
    \caption{ Solutions of the system (\ref{1}) for various values of $b$ inside and outside the stable region.}
    \label{pa1}
\end{figure}

\begin{ex}
    Let $\alpha=0.4$, $\beta=0.2$ and $a=1$. According to Fig. \ref{beta_0.2}, we take $b=0.5+i$ in the stable region and the initial condition as x(0)=0.1. The Fig. \ref{pa1a} shows that the solution of system (\ref{1}) converges to $0$. On the other hand, if we take $b=-1.5+i$ and $b=2.5+i$ outside the stable region, the system (\ref{1}) is unbounded ( cf. Fig.\ref{pa1b} and Fig.\ref{pa1c} ).
\end{ex}

\begin{figure}[H]
    \centering
    \subfigure[$b=-0.7-0.1i$]{
\includegraphics[height=2.5in,width=2.5in,keepaspectratio]{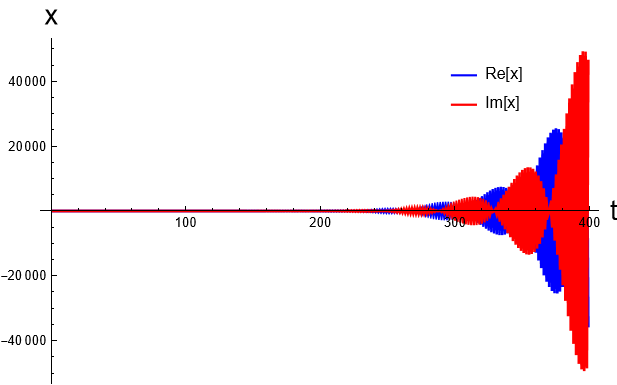}
        \label{pa2a}
    } \hspace{0.3cm}
    \subfigure[$b=1.2+0.3i$]{
\includegraphics[height=2.5in,width=2.5in,keepaspectratio]{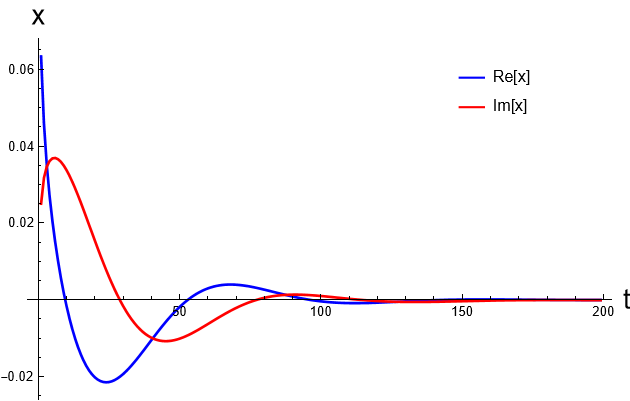}
        \label{pa2b}
    } \hspace{0.3cm}
    \subfigure[$b=1.4+0.1i$]{
\includegraphics[height=2.5in,width=2.5in,keepaspectratio]{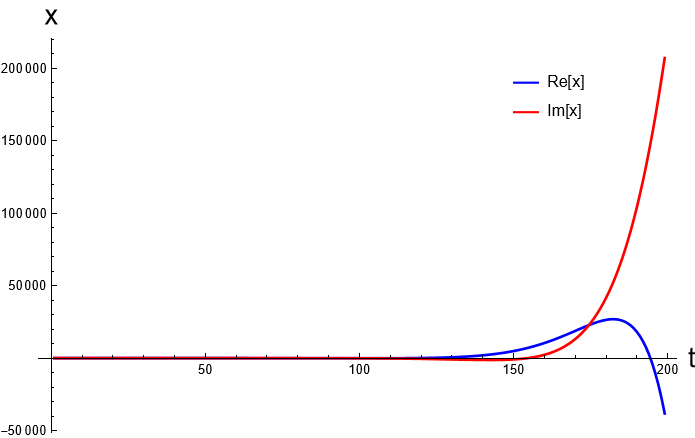}
        \label{pa2c}
    }
    \caption{ Solutions of the system (\ref{1}) for various values of $b$ inside and outside the stable region.}
    \label{pa2}
\end{figure}

\begin{ex}
    Consider $\alpha=0.5$, $\beta=0.3$ and $a=0.2$. From the Fig. \ref{alpha_0.5}, we take $b=1.2+0.3i$ in the stable region and the initial condition as x(0)=0.1. The solution of system (\ref{1}) converges to $0$ as shown in the Fig. \ref{pa2b}. However, if we take $b=-0.7-0.1i$ and $b=1.4+0.1i$ outside the stable region, the system (\ref{1}) is unbounded ( cf. Fig.\ref{pa2a} and Fig.\ref{pa2c} ).
\end{ex}

\subsection{Bifurcation analysis}

 For $a<0$, there are various bifurcations. We give the details of this rich dynamics in the following cases.\\

 \par \textbf{Case(I): Bifurcations for $0.5<\beta<\alpha<1$}\\

  \par
\begin{figure}[H]
    \centering
    \subfigure[$a_1$]{
\includegraphics[height=2.1in,width=2.1in,keepaspectratio]{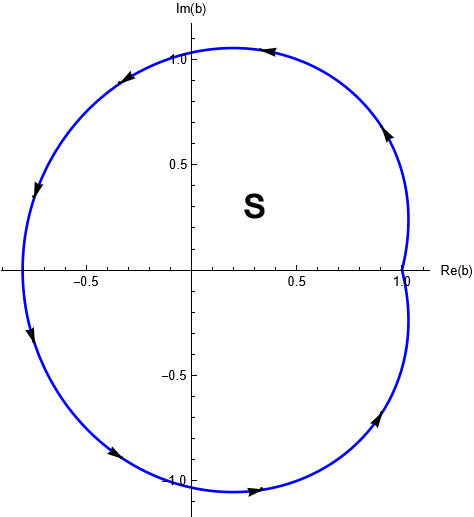}
        \label{case1a}
    } \hspace{0.3cm}
    \subfigure[$a_2<a<a_1$]{
\includegraphics[height=2.1in,width=3.1in,keepaspectratio]{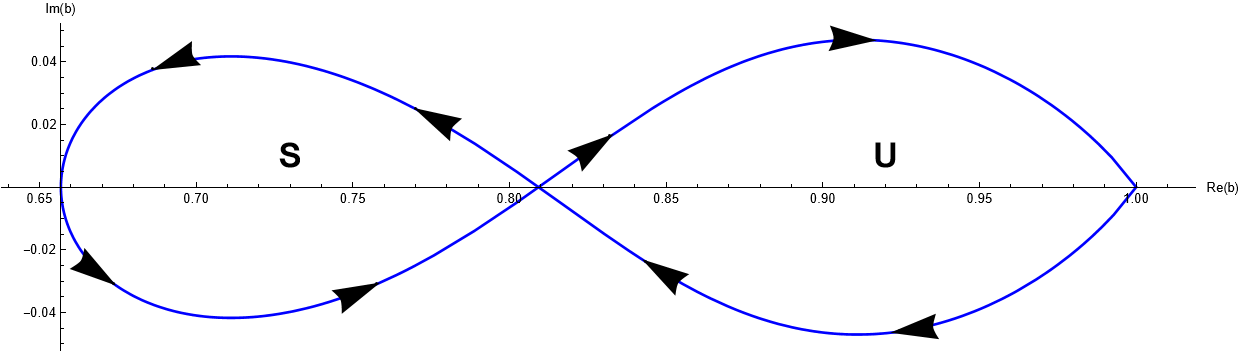}
        \label{case1b}
    } \hspace{0.3cm}
    \subfigure[$a_2$]{
\includegraphics[height=2.1in,width=2.1in,keepaspectratio]{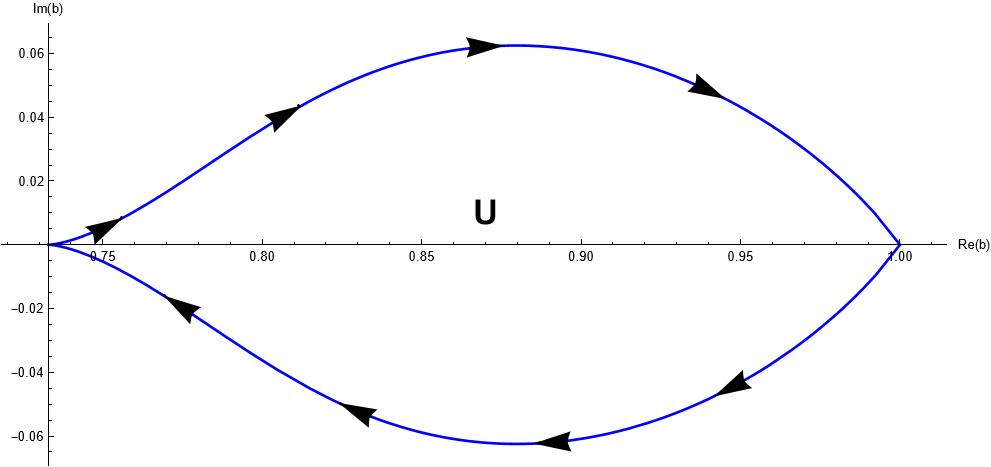}
        \label{case1c}
    } \hspace{0.3cm}
    \subfigure[$a_3<a<a_2$]{
\includegraphics[height=2.1in,width=2.1in,keepaspectratio]{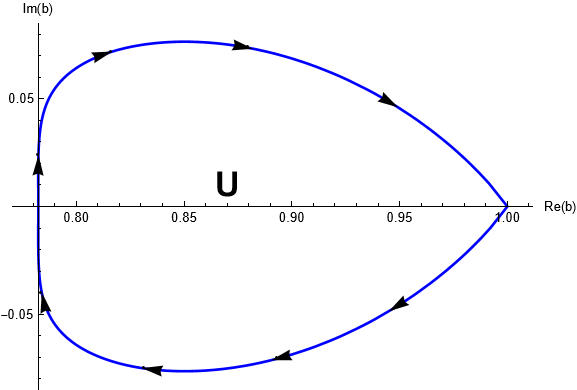}
        \label{case1d}
    } \hspace{0.3cm}
    \subfigure[$a_3$]{
\includegraphics[height=2.1in,width=2.1in,keepaspectratio]{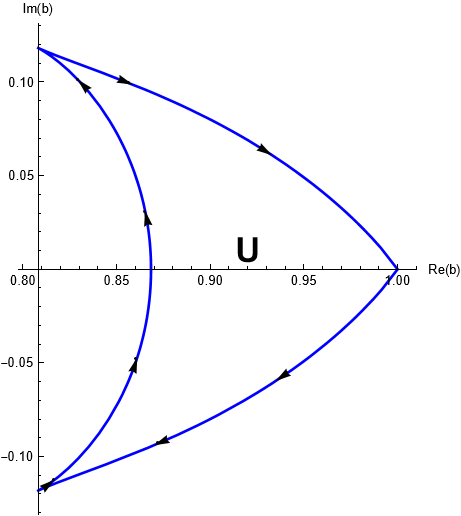}
        \label{case1e}
    } \hspace{0.6cm}
    \subfigure[$a_4<a<a_3$]{
\includegraphics[height=2.1in,width=2.1in,keepaspectratio]{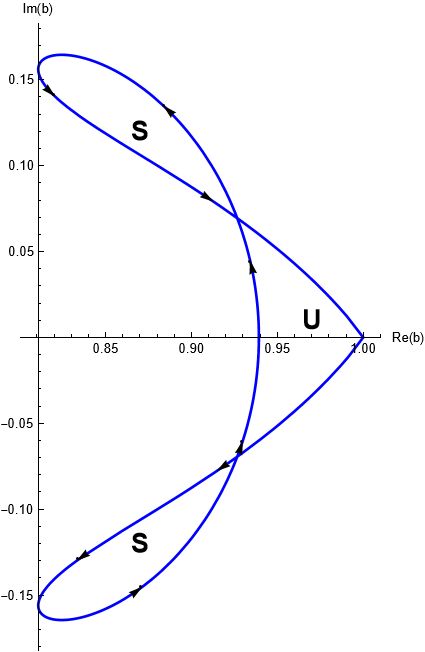}
        \label{case1f}
    } \hspace{0.6cm}
    \subfigure[$a_4$]{
\includegraphics[height=2.1in,width=2.1in,keepaspectratio]{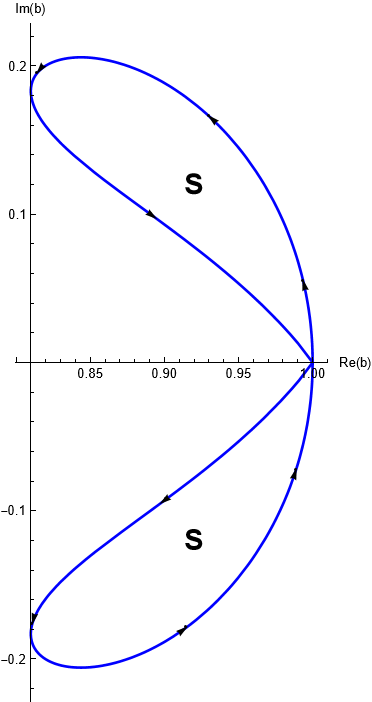}
        \label{case1g}
    } \hspace{0.3cm}
    \subfigure[$a<a_4$]{
\includegraphics[height=2.1in,width=2.1in,keepaspectratio]{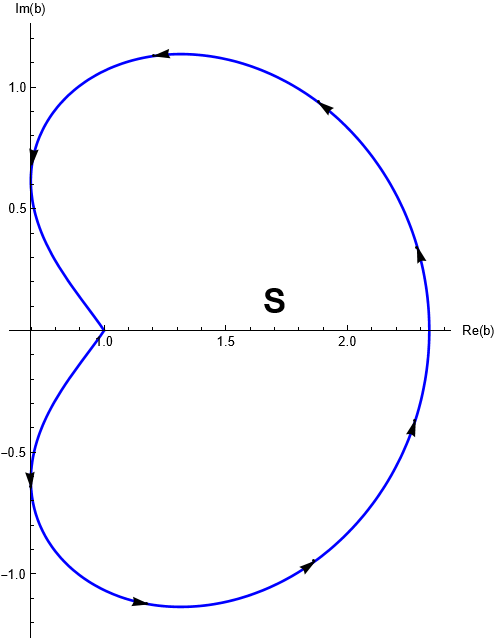}
        \label{case1h}
    }
    \caption{ (a)-(h) are representative stability diagrams for $0.5<\beta<\alpha<1$ and different values of $a$. In figures, S: stable, U: unstable.}
    \label{case1all}
\end{figure}

 It is observed that, as we decrease the value of $a$, the region bounded by $\gamma(\theta)$ divides into sub-regions. This division is due to self-intersections in $\gamma$. Such multiple points in $\gamma$ emerge due to cusps in it. The first bifurcation value occurs at $a_1=0$. At $a=a_1$, the stable region is a cardioid shape ( see Fig. \ref{case1a} ). For $a<a_1=0$ (sufficiently close to $a_1$), the curve $\gamma(\theta)$ becomes a multiple curve with a self-intersection ( see Fig. \ref{case1b} ). It gives rise to two closed regions. The left one is bounded by a positively oriented part of the curve $\gamma(\theta)$ and hence it is a stable region for the system (\ref{1}). The right part of the curve $\gamma(\theta)$ is unstable because it is bounded by the negatively oriented part of the curve $\gamma(\theta)$. Therefore, the unstable region starts at $a_1$.\\
 \par The curve $\gamma$ has a cusp if $\gamma^{'}(\theta)=0$.\\
 i.e., if
\begin{equation}
\begin{aligned}
&2^{-1 + \alpha} \alpha \cos\left(\frac{\theta}{2}\right) 
\cos\left(\frac{1}{2} \alpha (\pi - \theta) + \theta\right) 
\sin\left(\frac{\theta}{2}\right)^{-1 + \alpha} \\
&+ 2^{-1 + \beta} a \beta \cos\left(\frac{\theta}{2}\right) 
\cos\left(\frac{1}{2} \beta (\pi - \theta) + \theta\right) 
\sin\left(\frac{\theta}{2}\right)^{-1 + \beta} \\
&- 2^{\alpha} \left(1 - \frac{\alpha}{2}\right) 
\sin\left(\frac{\theta}{2}\right)^{\alpha} 
\sin\left(\frac{1}{2} \alpha (\pi - \theta) + \theta\right) \\
&- 2^{\beta} a \left(1 - \frac{\beta}{2}\right) 
\sin\left(\frac{\theta}{2}\right)^{\beta} 
\sin\left(\frac{1}{2} \beta (\pi - \theta) + \theta\right)=0
\end{aligned} \label{derre}
\end{equation}

and 
\begin{equation}
\begin{aligned}
&2^{\alpha} \left(1 - \frac{\alpha}{2}\right) 
\cos\left( \frac{1}{2} \alpha (\pi - \theta) + \theta \right) 
\sin\left( \frac{\theta}{2} \right)^{\alpha} \\
&+ 2^{\beta} a \left(1 - \frac{\beta}{2}\right) 
\cos\left( \frac{1}{2} \beta (\pi - \theta) + \theta \right) 
\sin\left( \frac{\theta}{2} \right)^{\beta} \\
&+ 2^{-1 + \alpha} \alpha \cos\left( \frac{\theta}{2} \right) 
\sin\left( \frac{\theta}{2} \right)^{-1 + \alpha} 
\sin\left( \frac{1}{2} \alpha (\pi - \theta) + \theta \right) \\
&+ 2^{-1 + \beta} a \beta \cos\left( \frac{\theta}{2} \right) 
\sin\left( \frac{\theta}{2} \right)^{-1 + \beta} 
\sin\left( \frac{1}{2} \beta (\pi - \theta) + \theta \right)=0.
\end{aligned} \label{derim}
\end{equation}
There is always a cusp at $\theta=0$. This cusp generates multiple points in $\gamma$ for $a<0$. Therefore, $a_1$ is a bifurcation point. If $\theta=\pi$, then equation (\ref{derre}) is satisfied identically. Therefore, solving equation (\ref{derim}) for $\theta=\pi$, we get a second bifurcation value $a_2=-\frac{2^{\alpha-\beta}(\alpha-2)}{\beta-2}$. At $a=a_2$, the left stable region shrinks and vanishes ( see Fig. \ref{case1c} ). For $a<a_2$ (sufficiently close to $a_2$), the right unstable region grows ( see Fig. \ref{case1d} ). In this case, we don't have any stable region and hence the system (\ref{1}) is unstable for any $b\in \mathbb{C}$.\\
 \par The third bifurcation value is $a_3$. At $a=a_3$, there are two cusps in $\gamma(\theta)$ ( see Fig. \ref{case1e} ). Those cusps can be obtained by solving equation (\ref{derim}) for $a$, and substituting this $a$ value in equation(\ref{derre}) to find the value of $\theta$.\\ Thus
$a_3=-\frac{ 2^{\alpha - \beta} \, \sin\left(\frac{\theta}{2}\right)^{\alpha - \beta} \left[ \sin\left( \frac{ -\theta(-3 + \alpha) + \pi \alpha }{2} \right) + (\alpha - 1) \sin\left( \frac{ \theta + \pi \alpha - \theta \alpha }{2} \right) \right] }{ \sin\left( \frac{ -\theta(-3 + \beta) + \pi \beta }{2} \right) + (\beta - 1) \sin\left( \frac{ \theta + \pi \beta - \theta \beta }{2} \right) }
$, where $\theta \in (0,\pi)$ can be found out by solving the following the equation

   \begin{eqnarray}
&&\left( -2 + \alpha + \beta - \alpha \beta \right) 
\sin\!\left[ \frac{ (\pi - \theta)(\alpha - \beta) }{2} \right] 
+ (\beta - 1) 
\sin\!\left[ \frac{ \theta(-2 + \alpha - \beta) + \pi(-\alpha + \beta) }{2} \right] \nonumber \\
&&\quad + (\alpha - 1) 
\sin\!\left[ \frac{ \theta(2 + \alpha - \beta) + \pi(-\alpha + \beta) }{2} \right]=0.\label{solthe}
\end{eqnarray}
We used the Mathematica command FindRoot to solve eq. (\ref{solthe}) ( see dataset1 available at the link \url{https://github.com/Janardhan3233/2term} ). For $a<a_3$ (sufficiently close to $a_3$), the curve $\gamma(\theta)$ has three regions ( see Fig. \ref{case1f} ). Here, two disjoint stable regions emerge.\\
\par The fourth bifurcation value is $a_4$. At $a=a_4$, the unstable region vanishes and the stable regions merge with each other ( see Fig. \ref{case1g} ). At $a=a_4$, we can observe that $\gamma(0)=\gamma(\pi)$, here $\gamma(0)=1$ and $\gamma(\pi)=1 - 2^{\alpha} - a \, 2^{\beta} $. Solving this, we get $a_4=-2^{\alpha-\beta}$. The curve $\gamma(\theta)$ will remain cardioid for $a<a_4$ ( see Fig. \ref{case1h} ). Here $b=1-2^{\alpha}-a \, 2^{\beta}$ and $b=1$ (cusp) are the right and left endpoints, respectively.\\
\subsubsection{Validation of analytic results  with Numerical Experiments}
In the following table, we provide some examples validating our theory for $a \leq 0$. We fix $\alpha=0.9$, $\beta=0.6$ and consider the initial condition as $x(0)=0.1$. If the value of $b$ lies inside the stable region, then we get a solution to system (\ref{1}) converging to $0$. 

\begin{table}[h!]
\centering
\renewcommand{\arraystretch}{1.5}
\begin{tabular}{|c|c|c|c|c|}
\hline
\textbf{a} & \textbf{Boundary Curve} & \textbf{b} & \textbf{Comment} & \textbf{Solution Curve} \\ 
\hline

{0} 
& {\includegraphics[height=3.5cm]{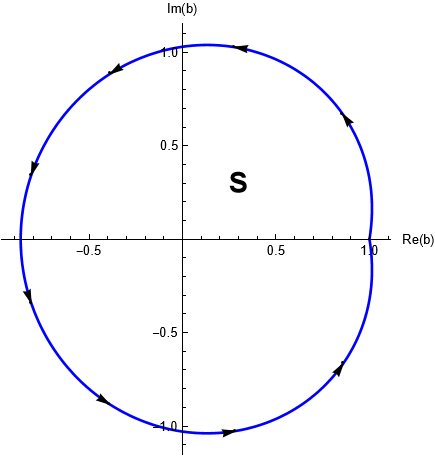}} 
& $0.5 + 0.5i$ & Stable & Fig. \ref{1c1a} \\ \cline{3-5}
& & $1-0.9i$ & Unstable & Fig. \ref{1c1b} \\ 
\hline

{-0.89}
& {\includegraphics[height=1.5cm]{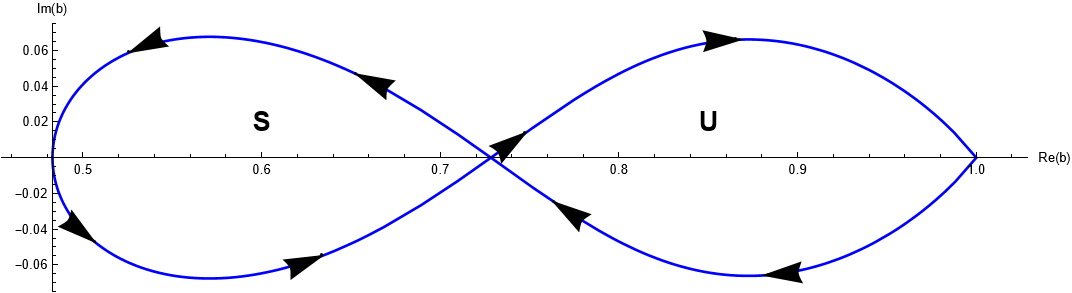}} 
& $0.6$ & Stable & Fig. \ref{1c1c} \\ \cline{3-5}
& & $0.9$ & Unstable & Fig. \ref{1c1d}\\ 
\hline

$-0.967328$ & \includegraphics[height=1.5cm]{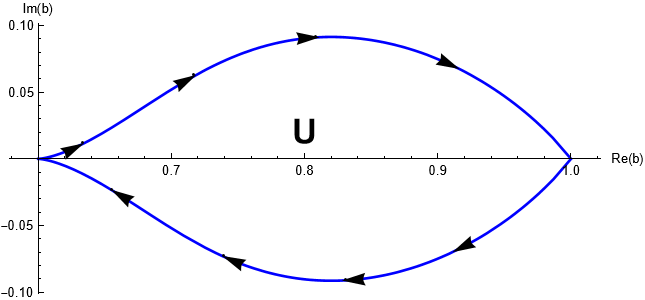} & $0.8$ & Unstable & Fig. \ref{1c1e} \\
\hline

{-1.17}
&{\includegraphics[height=4.5cm]{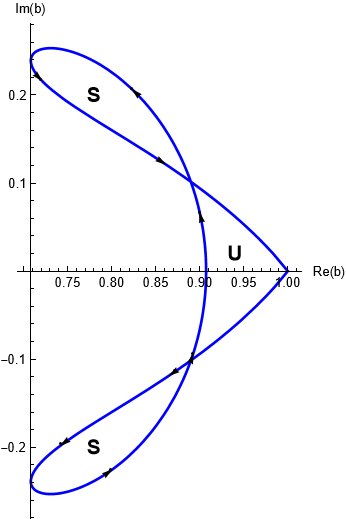}} 
& $0.8+0.2i$ & Stable & Fig. \ref{1c1f} \\ \cline{3-5}
& & $0.8-0.2i$ & Stable & Fig. \ref{1c1g}\\ \cline{3-5}
& & $0.95$ & Unstable & Fig. \ref{1c1h} \\ 
\hline

{-2.5} 
& {\includegraphics[height=3.5cm]{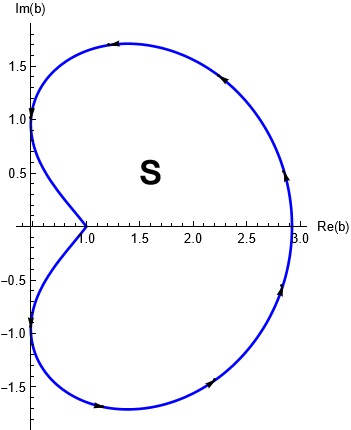}} 
& $1.5-i$ & Stable & Fig. \ref{1c1i}\\ \cline{3-5}
& & $0.5+0.2i$ & Unstable & Fig. \ref{1c1j} \\ 
\hline
\end{tabular}
\caption{Stability and solution curve data for different values of $a$ and $b$.}
\end{table}

\begin{figure}[H]
    \centering
    \subfigure[$a=0,\, b=0.5+0.5i$]{
\includegraphics[height=2.2in,width=2.2in,keepaspectratio]{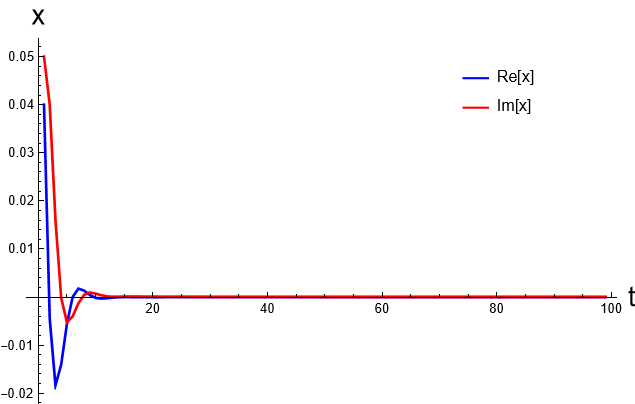}
        \label{1c1a}
    } \hspace{0.3cm}
    \subfigure[$a=0,\, b=1-0.9i$]{
\includegraphics[height=2.2in,width=2.2in,keepaspectratio]{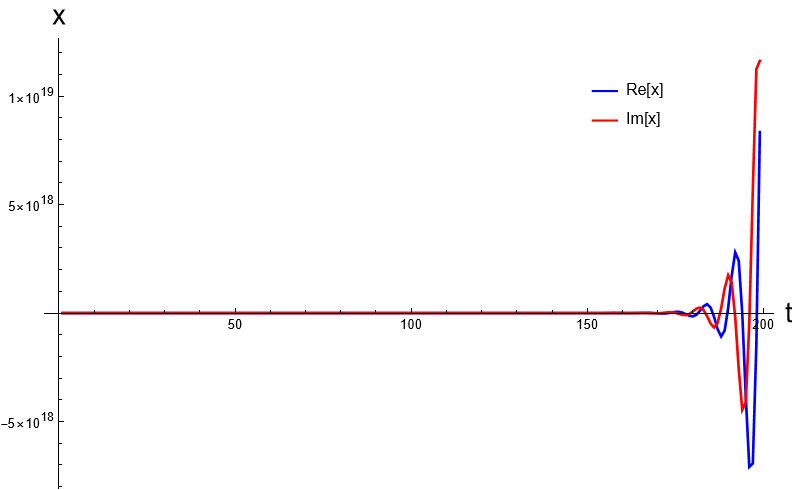}
        \label{1c1b}
    } \hspace{0.3cm}
    \subfigure[$a=-0.89, \, b=0.6$]{
\includegraphics[height=2.2in,width=2.2in,keepaspectratio]{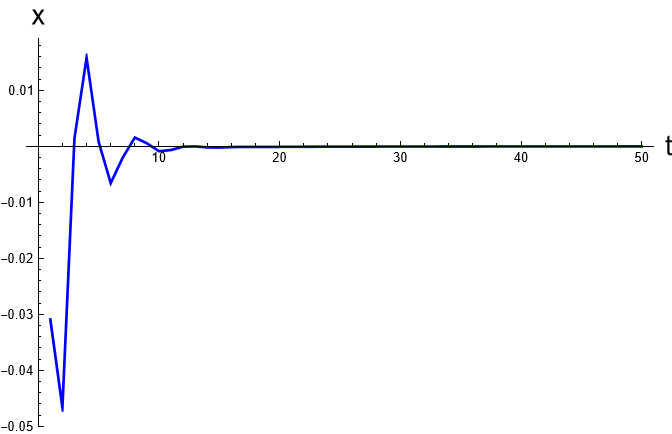}
        \label{1c1c}
    } \hspace{0.3cm}
    \subfigure[$a=-0.89, \, b=0.9$]{
\includegraphics[height=2.2in,width=2.2in,keepaspectratio]{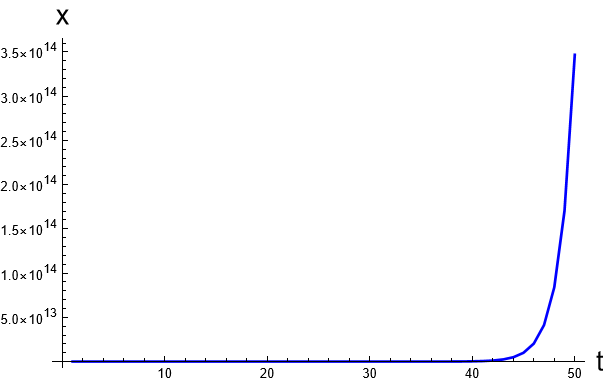}
        \label{1c1d}
    } \hspace{0.3cm}
    \subfigure[$a=-0.967328, \, b=0.8$]{
\includegraphics[height=2.2in,width=2.2in,keepaspectratio]{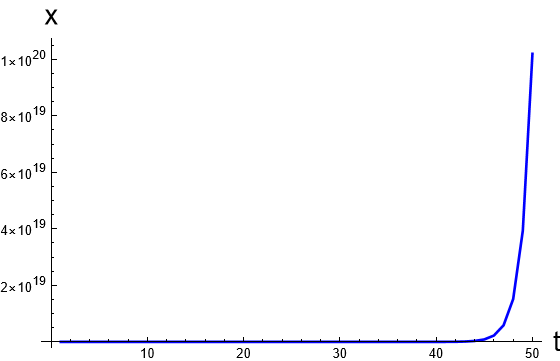}
        \label{1c1e}
    } \hspace{0.3cm}
    \subfigure[$a=-1.17, \, b=0.8+0.2i$]{
\includegraphics[height=2.2in,width=2.2in,keepaspectratio]{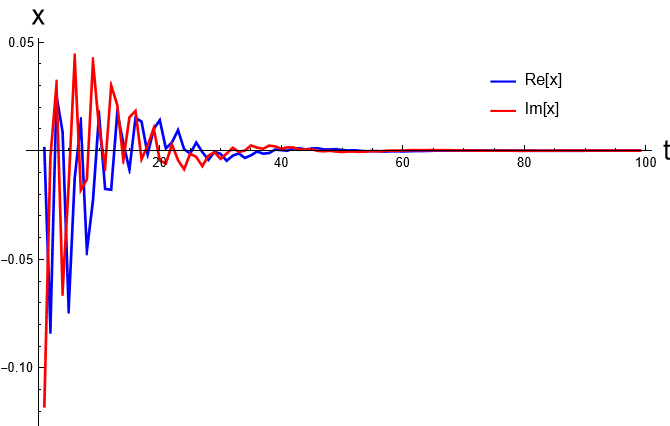}
        \label{1c1f}
    } \hspace{0.3cm}
    \subfigure[$a=-1.17, \, b=0.8-0.2i$]{
\includegraphics[height=2.2in,width=2.2in,keepaspectratio]{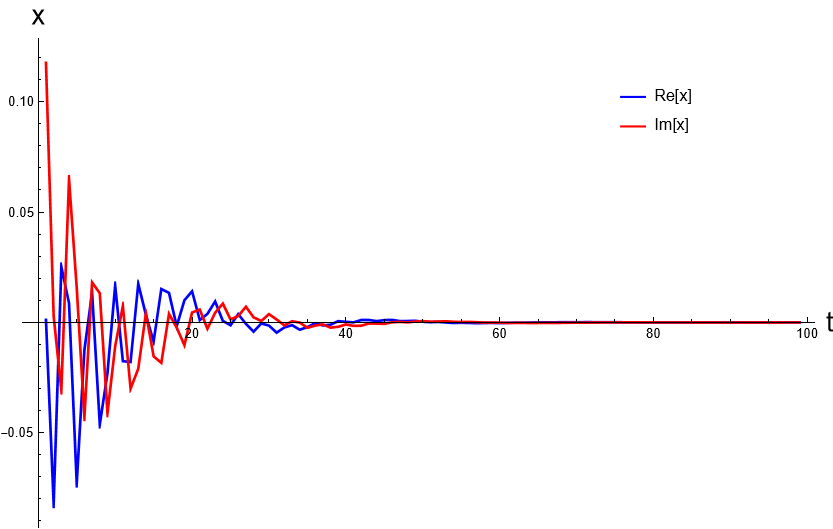}
        \label{1c1g}
    } \hspace{0.3cm}
    \subfigure[$a=-1.17, \, b=0.95$]{
\includegraphics[height=2.2in,width=2.2in,keepaspectratio]{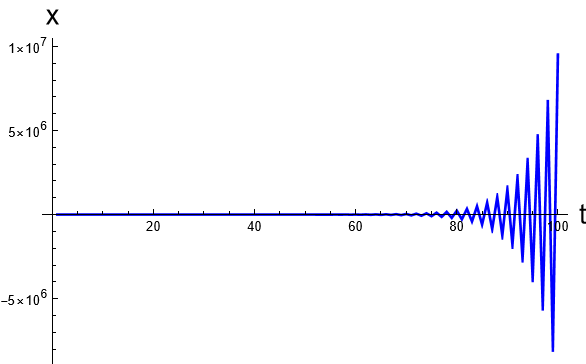}
        \label{1c1h}
    } \hspace{0.3cm}
    \subfigure[$a=-2.5, \, b=1.5-i$]{
\includegraphics[height=2.2in,width=2.2in,keepaspectratio]{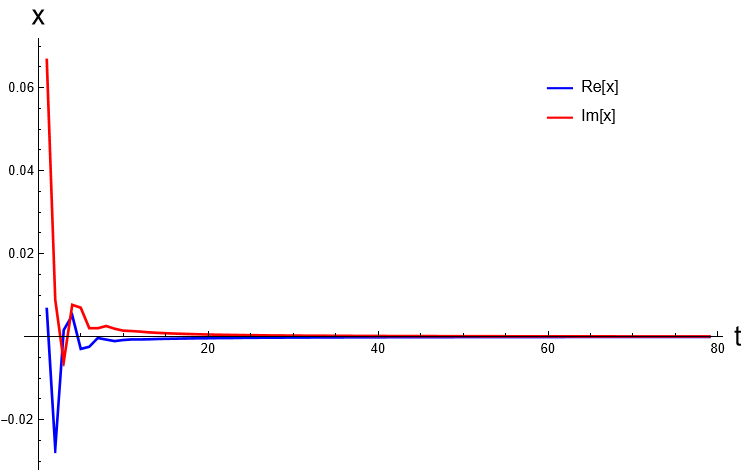}
        \label{1c1i}
    } \hspace{0.3cm}
    \subfigure[$a=-2.5, \, b=0.5+0.2i$]{
\includegraphics[height=2.2in,width=2.2in,keepaspectratio]{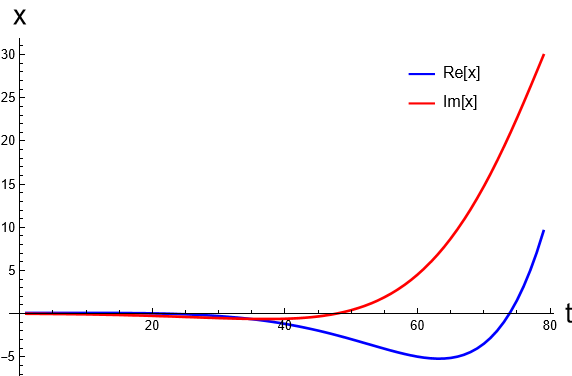}
        \label{1c1j}
    }
    \caption{ (a)-(j) are the solutions of the system (\ref{1}) for $\alpha=0.9$, $\beta=0.6$ and different values of $a$ and $b$.}
    \label{1c1}
\end{figure}

In the following, we discuss the case of stability with fixed $\alpha \, , \beta \in (0,0.5)$.\\ 

  \par \textbf{Case-II: Bifurcations for $0<\beta<\alpha<0.5$}\\

\begin{figure}[H]
\centering
\begin{tabular}{ccc} % 3 figures per row
\subfigure[$a_1$\label{case2a}]{\includegraphics[width=0.2\textwidth]{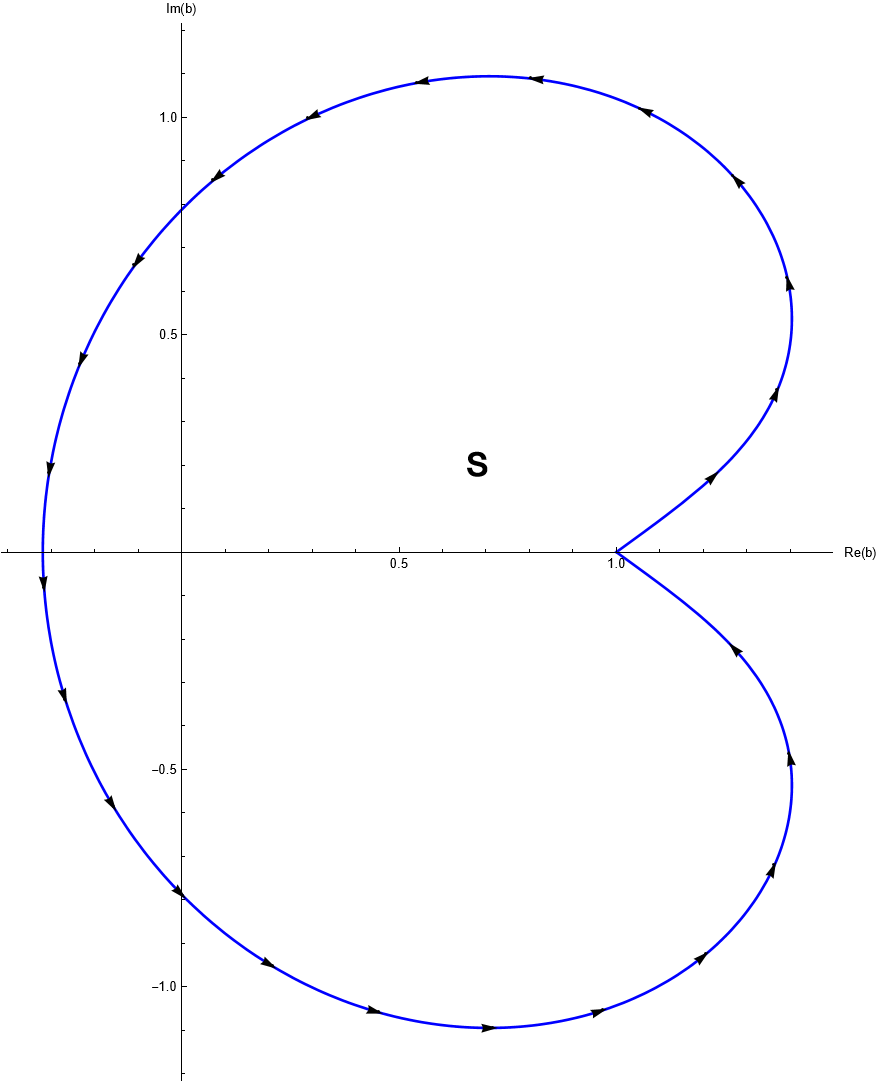}} &
\subfigure[$a_2<a<a_1$\label{case2b}]{\includegraphics[width=0.25\textwidth]{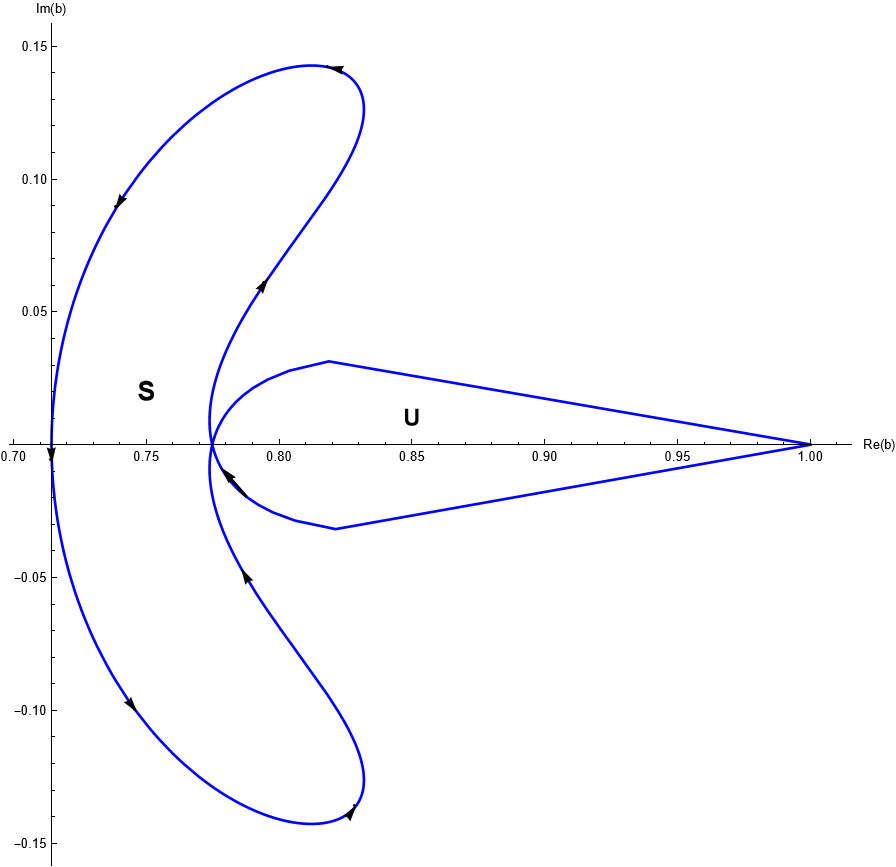}} &
\subfigure[$a_2$\label{case2c}]{\includegraphics[width=0.3\textwidth]{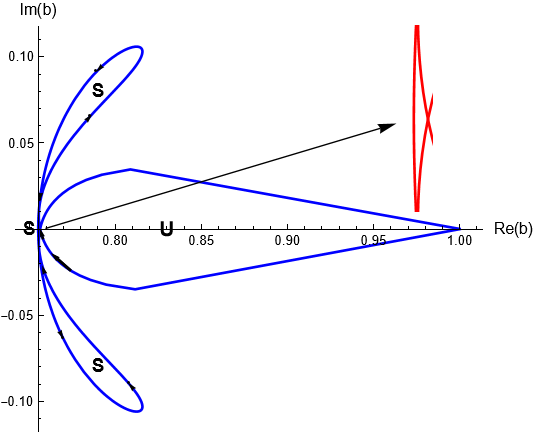}} \\

\subfigure[$a_3<a<a_2$\label{case2d}]{\includegraphics[width=0.25\textwidth]{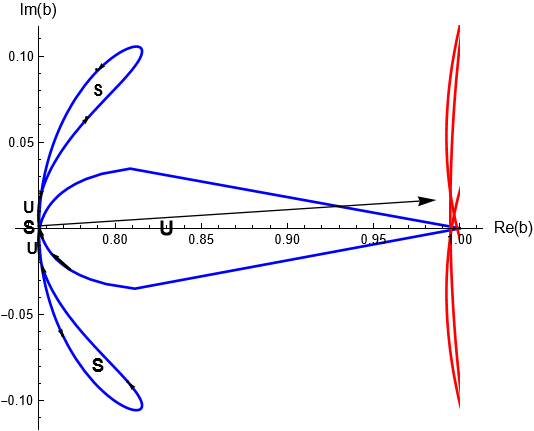}} &
\subfigure[$a_3$\label{case2e}]{\includegraphics[width=0.25\textwidth]{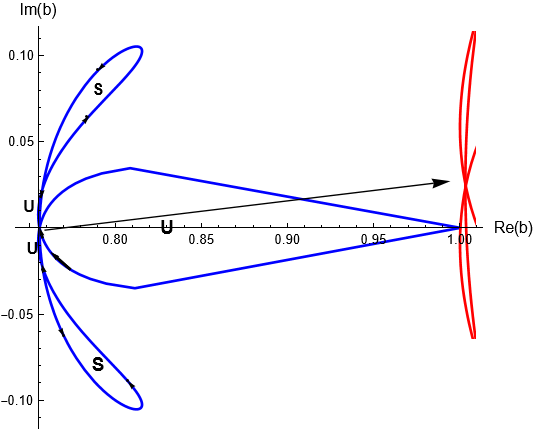}} &
\subfigure[$a_4<a<a_3$\label{case2f}]{\includegraphics[width=0.25\textwidth]{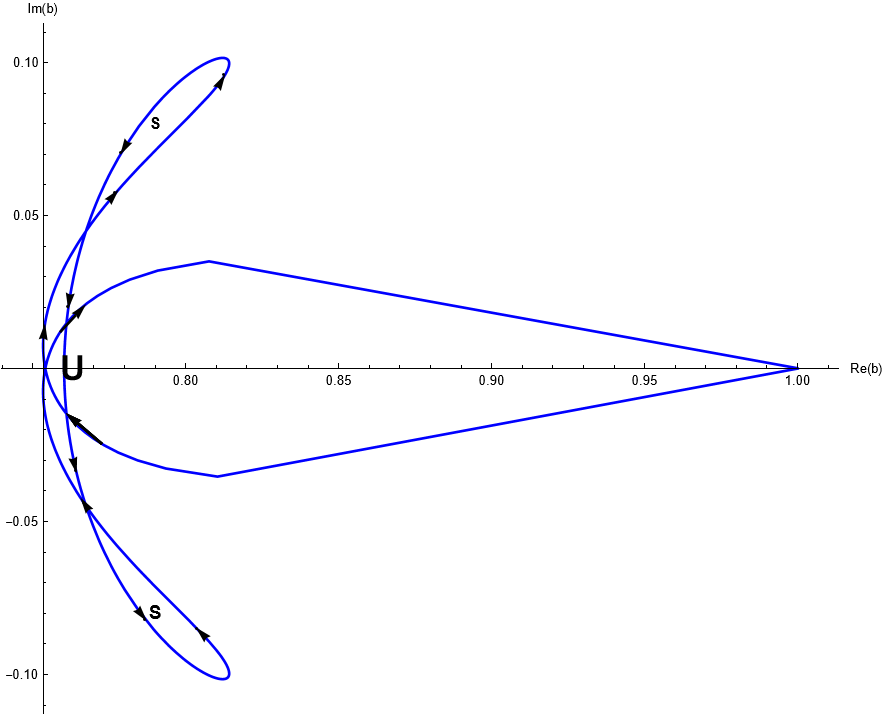}} \\

\subfigure[$a_4$\label{case2g}]{\includegraphics[width=0.25\textwidth]{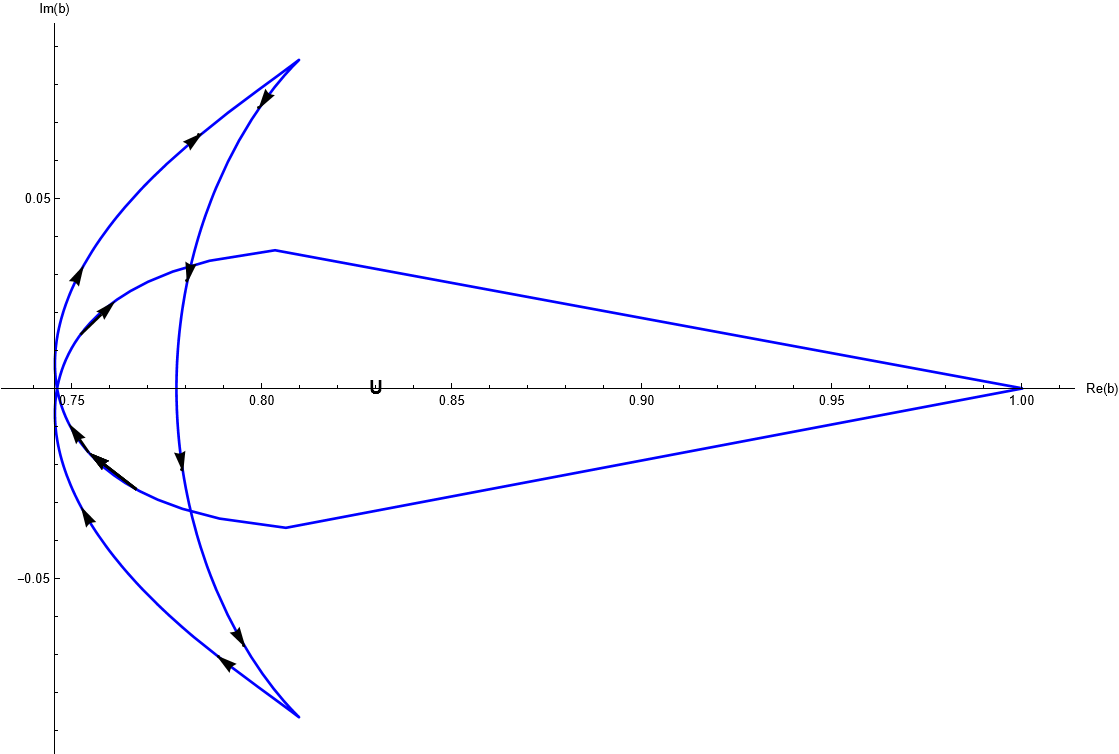}} &
\subfigure[$a_5<a<a_4$\label{case2h}]{\includegraphics[width=0.35\textwidth]{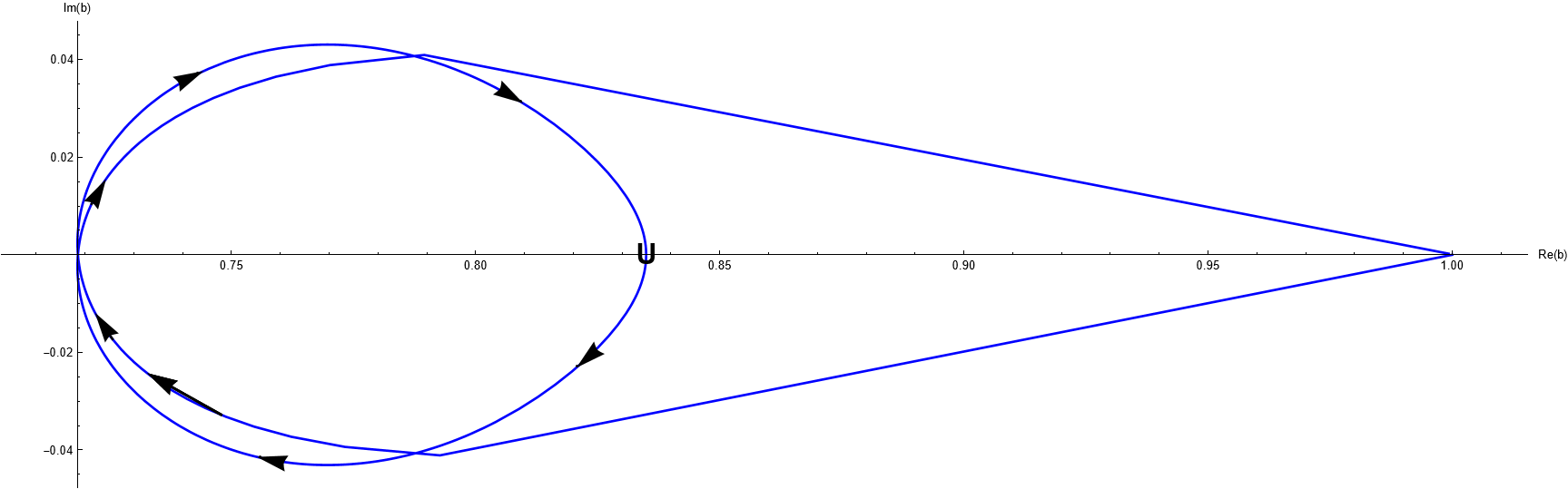}} &
\subfigure[$a_5$\label{case2i}]{\includegraphics[width=0.35\textwidth]{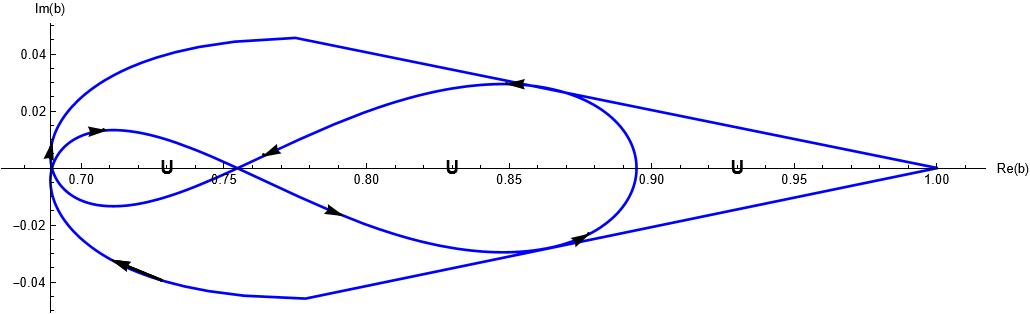}} \\

\subfigure[$a_6<a<a_5$\label{case2j}]{\includegraphics[width=0.35\textwidth]{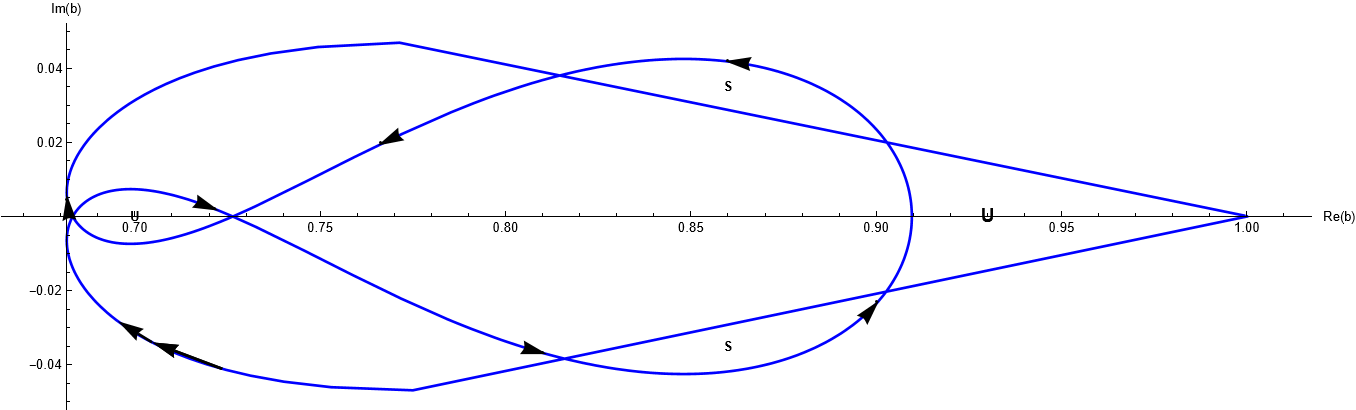}} &
\subfigure[$a_6$\label{case2k}]{\includegraphics[width=0.35\textwidth]{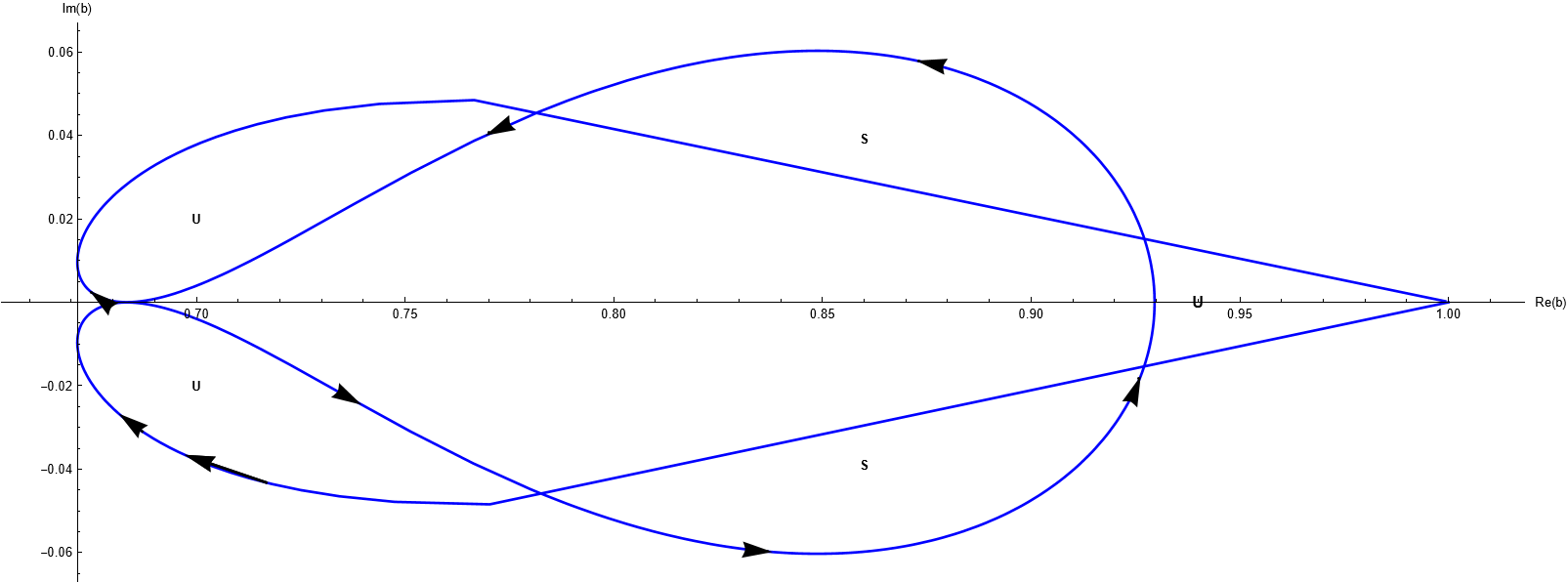}} &
\subfigure[$a_7<a<a_6$\label{case2l}]{\includegraphics[width=0.35\textwidth]{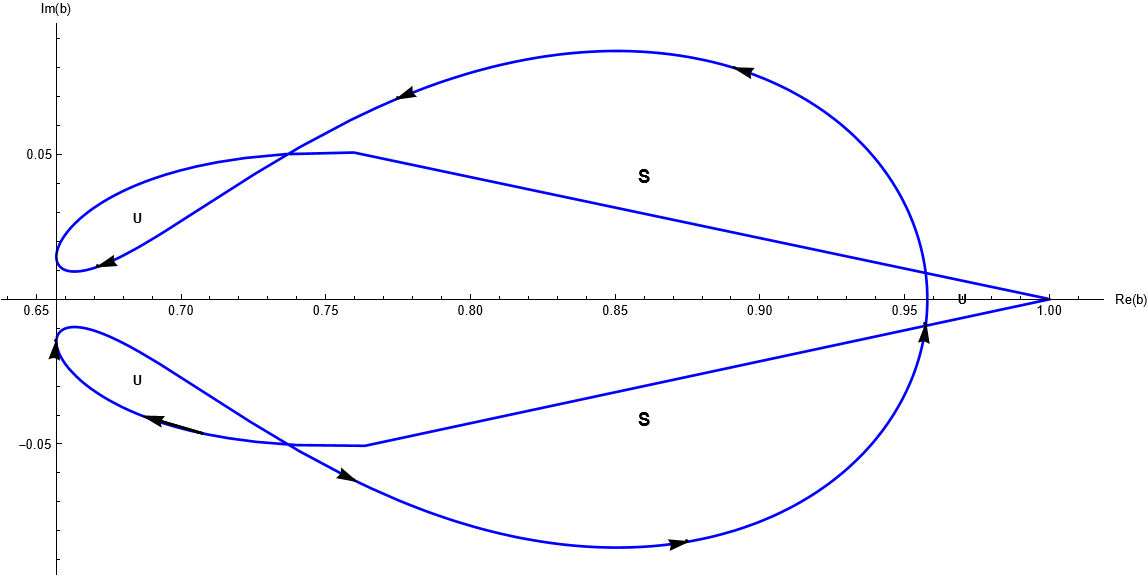}} \\

\subfigure[$a_7$\label{case2m}]{\includegraphics[width=0.3\textwidth]{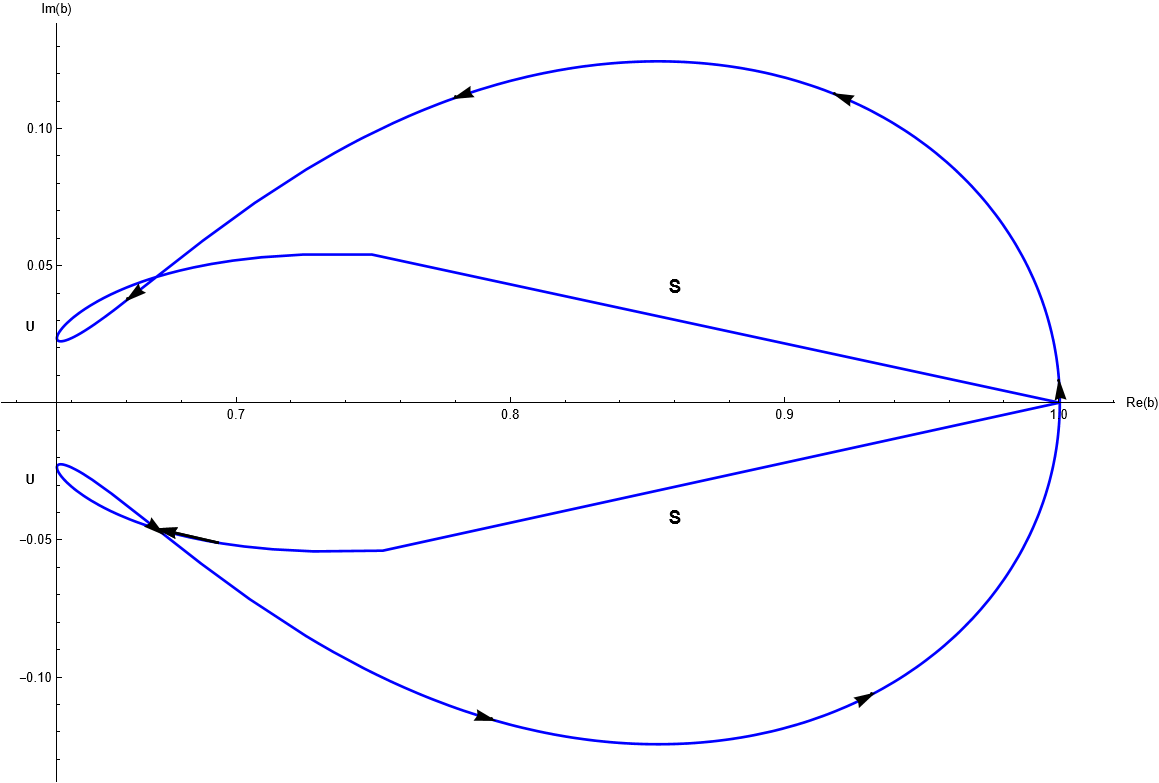}} &
\subfigure[$a_8<a<a_7$\label{case2n}]{\includegraphics[width=0.25\textwidth]{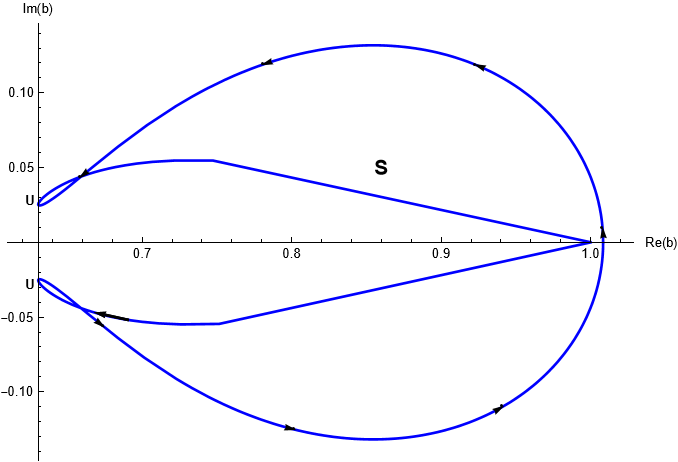}} &
\subfigure[$a_8$\label{case2o}]{\includegraphics[width=0.2\textwidth]{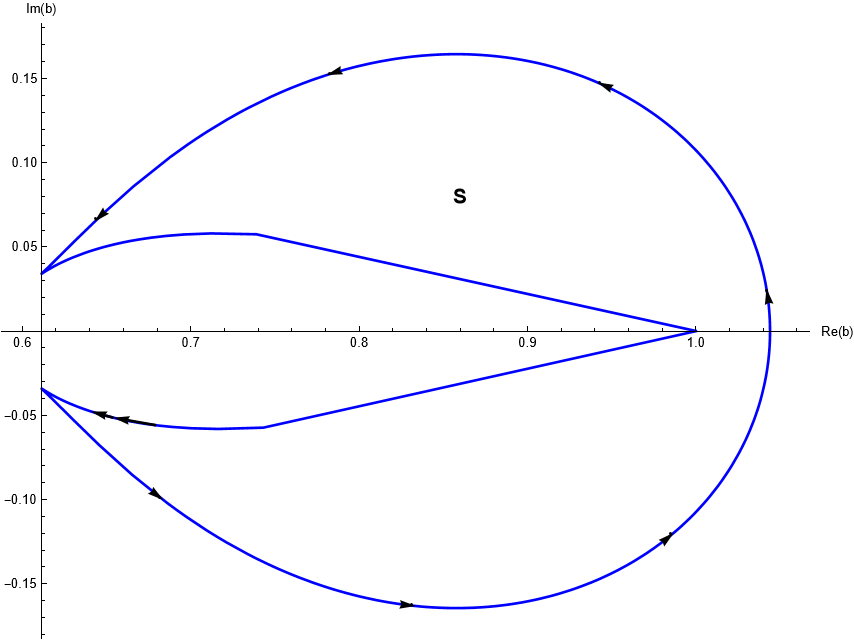}} \\

\subfigure[$a<a_8$\label{case2p}]{\includegraphics[width=0.2\textwidth]{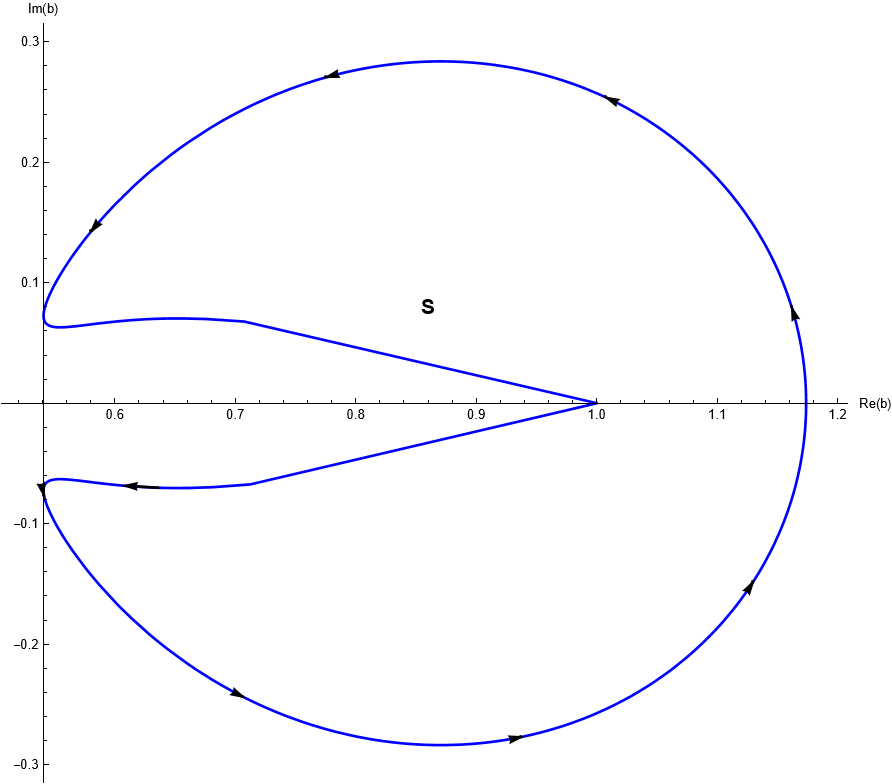}} 

\end{tabular}
\caption{Stability diagrams for $\alpha=0.4$, $\beta=0.2$ and different values of $a$.}
\label{case2all}
\end{figure}

  For this case, we illustrate the bifurcation scenario by taking particular values $\alpha=0.4$ and $\beta=0.2$. We sketch the boundary curve $\gamma(\theta)$ for different values of parameter $a$ and provide the stable and unstable regions of the system (\ref{1}) for different values of $a$ ( cf. Fig. \ref{case2all} ). Here, the dynamics have the transitions from a single stable region of the curve $\gamma(\theta)$ to two sub-regions bifurcated by $a=a_1=0$. For $a<a_1$, say $a=-0.9$ ( Fig. \ref{case2b} ), $\gamma(\theta)$ is divided into one stable and one unstable region due to self-intersections of $\gamma$. Therefore, the unstable region starts at $a_1$. The next bifurcation occurs at $a_2=-0.93569$, where the previous single stable sub-region of $\gamma(\theta)$ is divided into three sub-stable regions ( Fig. \ref{case2c} ). After crossing the bifurcation value $a_2$, say $a=-0.9361$ ( Fig. \ref{case2d} ), the two new unstable sub-regions emerge. The third bifurcation value is $a_3=-0.93629$. At $a_3$, one of the three stable sub-regions vanishes ( Fig. \ref{case2e} ). For $a<a_3$, say $a=-0.94$ ( Fig. \ref{case2f} ), the unstable sub-regions merge with each other. The next bifurcation value is $a_4=-0.954579$. At $a_4$, the two unstable sub-regions vanish and produce two cusps in $\gamma(\theta)$, and everything is unstable for any complex value of $b$ ( Fig. \ref{case2g} ). After crossing $a_4$, say $a=-1.005$ ( Fig. \ref{case2h} ), only the unstable region remains, and it grows up to the next bifurcation value $a_5=-1.057$. At $a_5$, one more unstable sub-region emerges on the right ( Fig. \ref{case2i} ). For $a<a_5$, say $a=-1.07$ ( Fig. \ref{case2j} ), the two new distinct stable sub-regions emerge. The next bifurcation value is $a_6=-1.0874$. At $a_6$, the left unstable sub-region converts to two distinct unstable sub-regions ( Fig. \ref{case2k} ). For $a<a_6$, say $a=-1.11$ ( Fig. \ref{case2l} ), the stable sub-regions are growing and unstable regions are shrinking. The next bifurcation value occurs at $a_7=-1.1487$, here $\gamma(0)=\gamma(\pi)$. At $a_7$, one of the three unstable sub-regions vanishes on the right ( Fig. \ref{case2m} ). For $a<a_7$, say $a=-1.156$  ( Fig. \ref{case2n} ), the two stable sub-regions merge with each other. At the next bifurcation value $a_8=-1.187$, the unstable sub-regions shrink and vanish on the left( Fig. \ref{case2o} ). After $a_8$ say $a=-1.3$, we have a single stable region ( Fig. \ref{case2p} ). We provide some examples to support the above stability regions.

 We verified all these regions in the following table. We set $\alpha=0.4$, $\beta=0.2$.
\begin{table}[p]
\centering
\renewcommand{\arraystretch}{1.5}
\begin{tabular}{|c|c|c|c|c|}
\hline
\textbf{a} & \textbf{Ref. Boundary Curve} & \textbf{b} & \textbf{Stability Property } & \textbf{Fig. Reference} \\ 
\hline

$0$ & Fig. \ref{case2a} & $1.215+0.6561i$ & Stable & Fig. \ref{acase2a} \\  
\cline{3-5} 
& & $i$ & Unstable & Fig. \ref{acase2b} \\ 
\hline

$-0.9$& Fig. \ref{case2b} & $0.8+0.1i$ & Stable & Fig. \ref{bcase2a} \\ 
\cline{3-5}
& & $0.85$ & Unstable & Fig. \ref{bcase2b}\\ 
\hline

$-0.93569$ & Fig. \ref{case2c} & $0.8+0.09i$ & Stable & Fig. \ref{ccase2a} \\
\cline{3-5}
& & $0.756$ & Stable & Fig. \ref{ccase2b}\\
\cline{3-5}
& & $0.78-0.07i$ & Stable & Fig. \ref{ccase2c}\\
\cline{3-5}
& & $0.82+0.02i$ & Unstable & Fig. \ref{ccase2d}\\
\hline

$-0.9361$ & Fig. \ref{case2d} & $0.79+0.08i$ & Stable & Fig. \ref{dcase2a} \\ 
\cline{3-5}
& & $0.7559+0.0002775i$ & Stable & Fig. \ref{dcase2b}\\
\cline{3-5}
& & $0.7984-0.09044i$ & Stable & Fig. \ref{dcase2c} \\
\cline{3-5}
& & $0.7558+0.01087i$ & Unstable & Fig. \ref{dcase2d} \\ 
\cline{3-5}
& & $0.7556-0.007171i$ & Unstable & Fig. \ref{dcase2e} \\ 
\cline{3-5}
& & $0.8223-0.01474i$ & Unstable & Fig. \ref{dcase2f} \\ 
\hline

$-0.94$ & Fig. \ref{case2f} & $0.7967+0.08422i$ & Stable & Fig. \ref{fcase2a}\\ \cline{3-5}
& & $0.8058-0.0958i$ & Stable & Fig. \ref{fcase2b} \\ 
\cline{3-5}
& & $0.8$ & Unstable & Fig. \ref{fcase2c} \\
\hline

$-1.07$ & Fig. \ref{case2j} & $0.8672+0.03309i$ & Stable & Fig. \ref{jcase2a} \\
\cline{3-5}
& & $0.8464-0.03691i$ & Stable & Fig. \ref{jcase2b}\\
\cline{3-5}
& & $0.7$ & Unstable & Fig. \ref{jcase2c}\\
\cline{3-5}
& & $0.95$ & Unstable & Fig. \ref{jcase2d}\\
\hline

$-1.0874$ & Fig. \ref{case2k} & $0.8667+0.04332i$ & Stable & Fig. \ref{kcase2a} \\
\cline{3-5}
& & $0.8424-0.04686i$ & Stable & Fig. \ref{kcase2b}\\
\cline{3-5}
& & $0.7+0.02i$ & Unstable & Fig. \ref{kcase2c}\\
\cline{3-5}
& & $0.75-0.04i$ & Unstable & Fig. \ref{kcase2d}\\
\cline{3-5}
& & $0.95$ & Unstable & Fig. \ref{kcase2e}\\
\hline

$-1.1487$ & Fig. \ref{case2m} & $0.8325+0.07685i$ & Stable & Fig. \ref{mcase2a} \\
\cline{3-5}
& & $0.8616-0.0896i$ & Stable & Fig. \ref{mcase2b}\\
\cline{3-5}
& & $0.6458+0.03003i$ & Unstable & Fig. \ref{mcase2c}\\
\cline{3-5}
& & $0.6497-0.03331i$ & Unstable & Fig. \ref{mcase2d}\\
\hline
\end{tabular}
\caption{Stability and solution curve data for $\alpha=0.4$, $\beta=0.2$ and different values of $a$ and $b$.}
\end{table}

\begin{figure}[p]
    \centering
    \subfigure[$a=0, b=1.215+0.6561i$]{
\includegraphics[height=1.5in,width=1.5in,keepaspectratio]{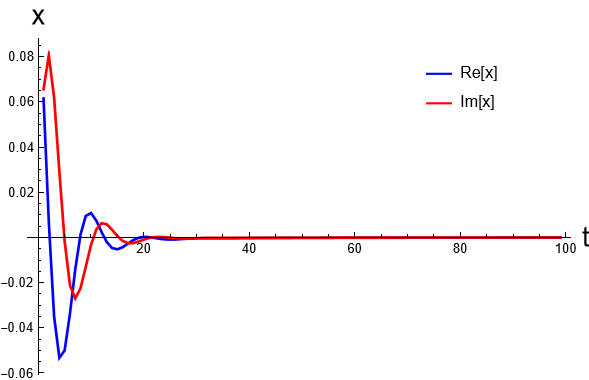}
        \label{acase2a}
    } \hspace{0.3cm}
    \subfigure[$a=0,\, b=i$]{
\includegraphics[height=1.5in,width=1.5in,keepaspectratio]{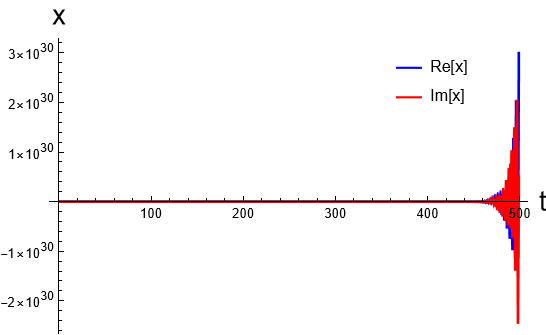}
        \label{acase2b}
    } \hspace{0.3cm}
    \subfigure[$a=-0.9, \, b=0.8+0.1i$]{
\includegraphics[height=1.5in,width=1.5in,keepaspectratio]{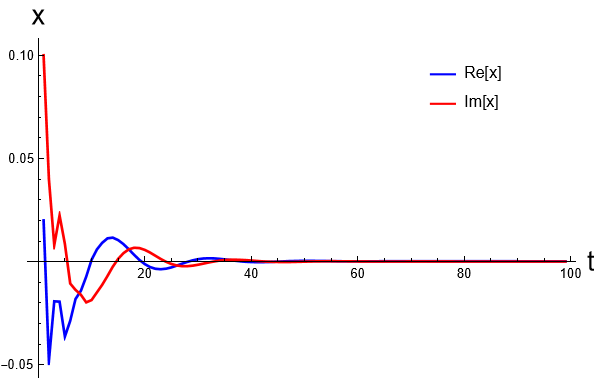}
        \label{bcase2a}
    } \hspace{0.3cm}
    \subfigure[$a=-0.9, \, b=0.85$]{
\includegraphics[height=1.5in,width=1.5in,keepaspectratio]{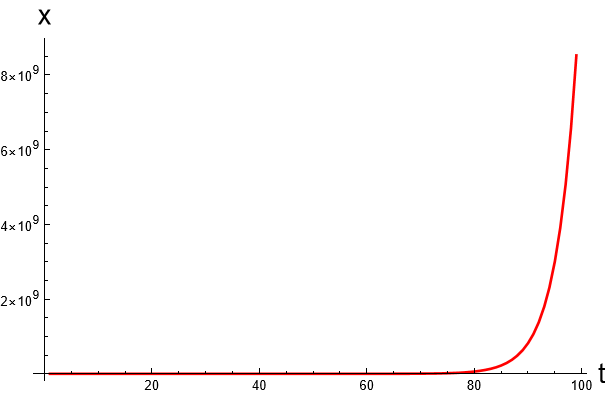}
        \label{bcase2b}
    } \hspace{0.3cm}
    \subfigure[$a=-0.93569, b=0.8+0.09i$]{
\includegraphics[height=1.5in,width=1.5in,keepaspectratio]{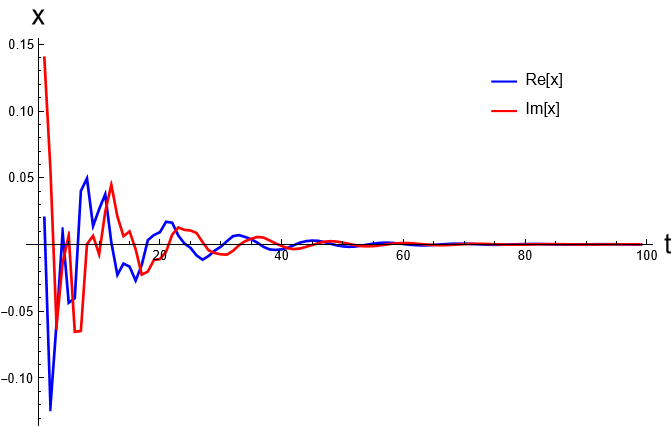}
        \label{ccase2a}
    } \hspace{0.3cm}
    \subfigure[$a=-0.93569, \, b=0.756$]{
\includegraphics[height=1.5in,width=1.5in,keepaspectratio]{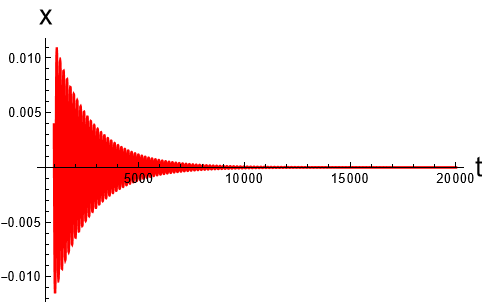}
        \label{ccase2b}
    } \hspace{0.3cm}
    \subfigure[$a=-0.93569, b=0.78-0.07i$]{
\includegraphics[height=1.5in,width=1.5in,keepaspectratio]{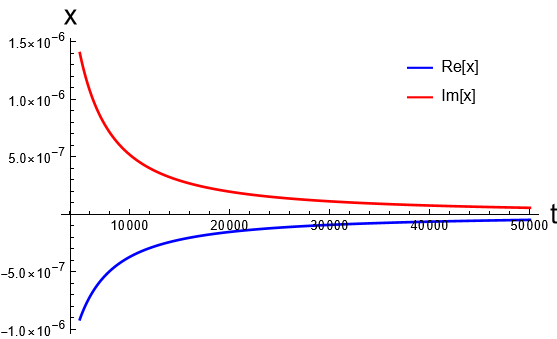}
        \label{ccase2c}
    } \hspace{0.3cm}
    \subfigure[$a=-0.93569, b=0.82+0.02i$]{
\includegraphics[height=1.5in,width=1.5in,keepaspectratio]{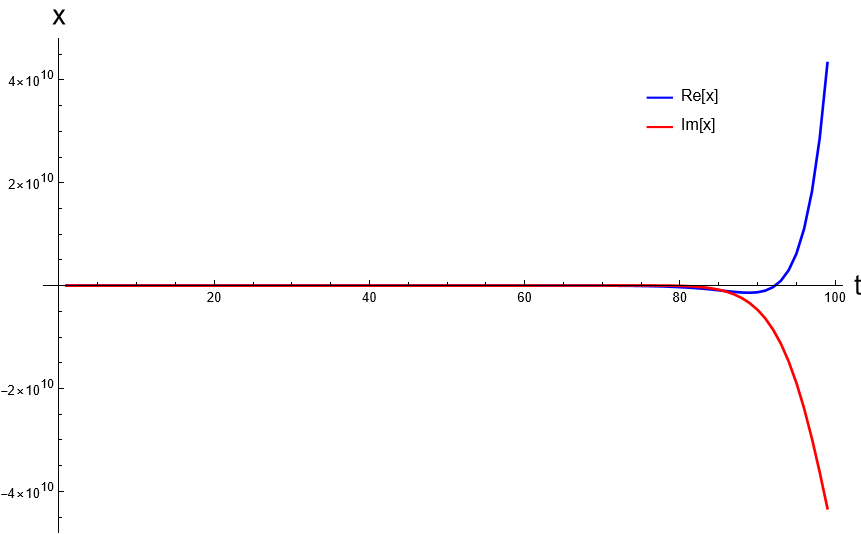}
        \label{ccase2d}
    } \hspace{0.3cm}
    \subfigure[$a=-0.9361, b=0.79+0.08i$]{
\includegraphics[height=1.5in,width=1.5in,keepaspectratio]{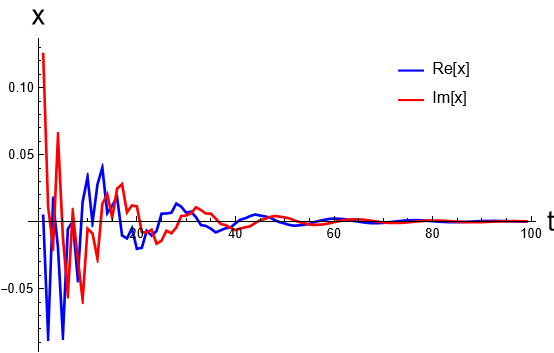}
        \label{dcase2a}
    } \hspace{0.3cm}
    \subfigure[$a=-0.9361, b=0.7559+0.0002775i$]{
\includegraphics[height=1.5in,width=1.5in,keepaspectratio]{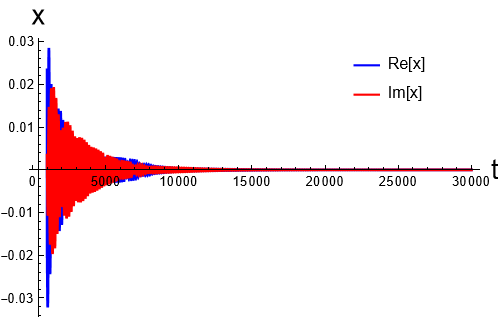}
        \label{dcase2b}
    } \hspace{0.3cm}
    \subfigure[$a=-0.9361, b=0.7984-0.09044i$]{
\includegraphics[height=1.5in,width=1.5in,keepaspectratio]{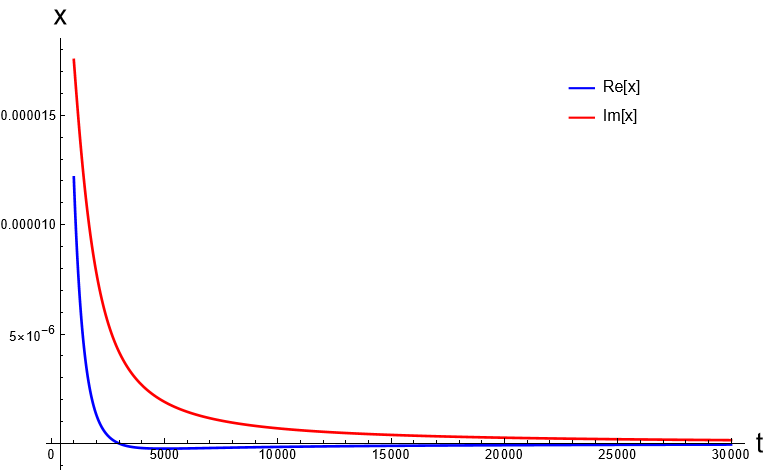}
        \label{dcase2c}
    } \hspace{0.3cm}
    \subfigure[$a=-0.9361, b=0.75580+0.01087i$]{
\includegraphics[height=1.5in,width=1.5in,keepaspectratio]{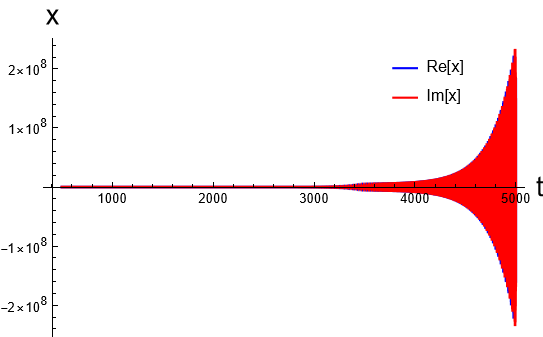}
        \label{dcase2d}
    } \hspace{0.3cm}
    \subfigure[$a=-0.9361, b=0.7556-0.007171i$]{
\includegraphics[height=1.5in,width=1.5in,keepaspectratio]{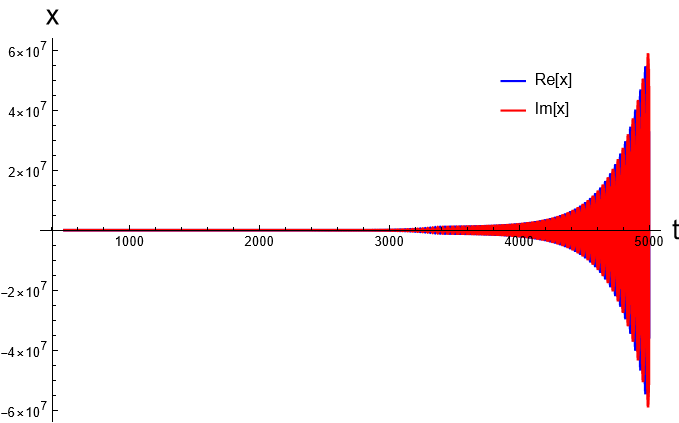}
        \label{dcase2e}
    } \hspace{0.3cm}
    \subfigure[$a=-0.9361, b=0.8223-0.01474i$]{
\includegraphics[height=1.5in,width=1.5in,keepaspectratio]{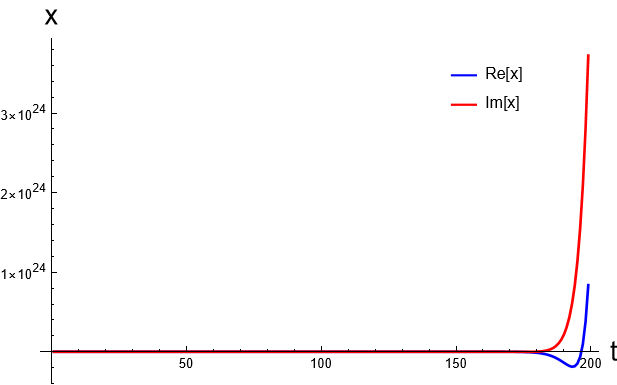}
        \label{dcase2f}
    } \hspace{0.3cm}
    \subfigure[$a=-0.94, b=0.7967+0.08422i$]{
\includegraphics[height=1.5in,width=1.5in,keepaspectratio]{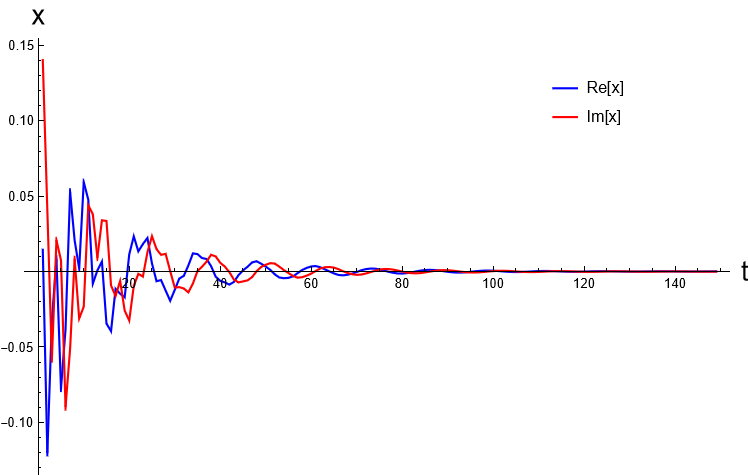}
        \label{fcase2a}
    } \hspace{0.3cm}
    \subfigure[$a=-0.94, b=0.8058-0.0958i$]{
\includegraphics[height=1.5in,width=1.5in,keepaspectratio]{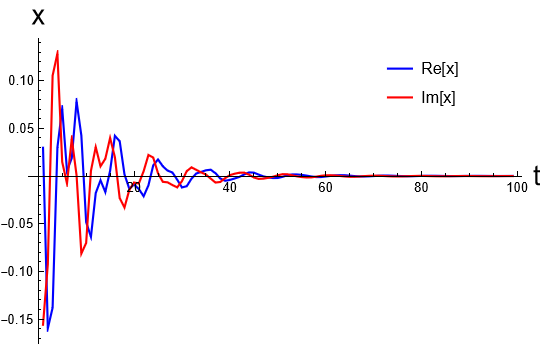}
        \label{fcase2b}
    } \hspace{0.3cm}
    \subfigure[$a=-0.94, \, b=0.8$]{
\includegraphics[height=1.5in,width=1.5in,keepaspectratio]{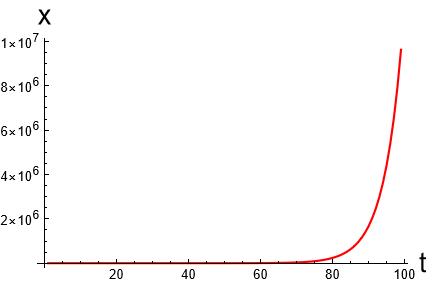}
        \label{fcase2c}
    } \hspace{0.3cm}
    \subfigure[$a=-1.07, b=0.8672+0.03309i$]{
\includegraphics[height=1.5in,width=1.5in,keepaspectratio]{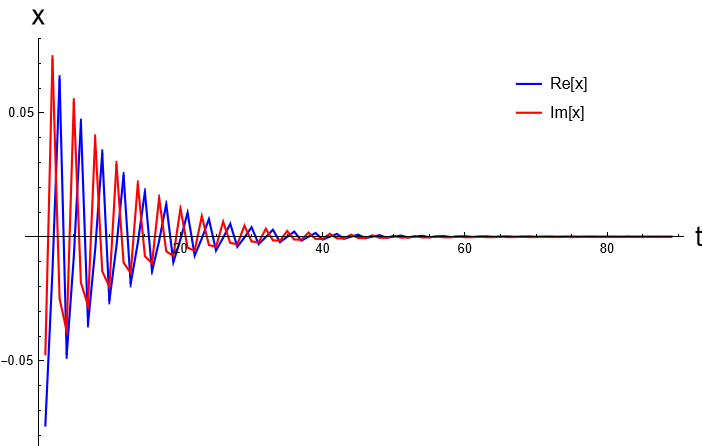}
        \label{jcase2a}
    
    }
    \caption{ Behavior of the solutions of the system (\ref{1}) for $\alpha=0.4$, $\beta=0.2$ and different values of $a$ and $b$.}
    \label{1c1}
\end{figure}

\begin{figure}[H]
    \centering
    \subfigure[$a=-1.07, b=0.8464-0.03691i$]{
\includegraphics[height=1.5in,width=1.5in,keepaspectratio]{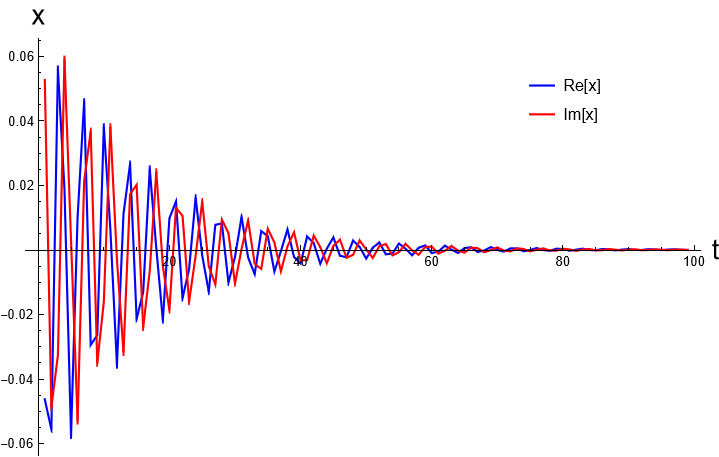}
        \label{jcase2b}
    } \hspace{0.3cm}
    \subfigure[$a=-1.07, \, b=0.7$]{
\includegraphics[height=1.5in,width=1.5in,keepaspectratio]{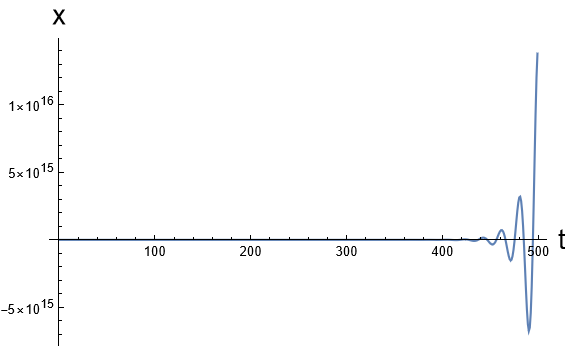}
        \label{jcase2c}
    } \hspace{0.3cm}
    \subfigure[$a=-1.07, \, b=0.95$]{
\includegraphics[height=1.5in,width=1.5in,keepaspectratio]{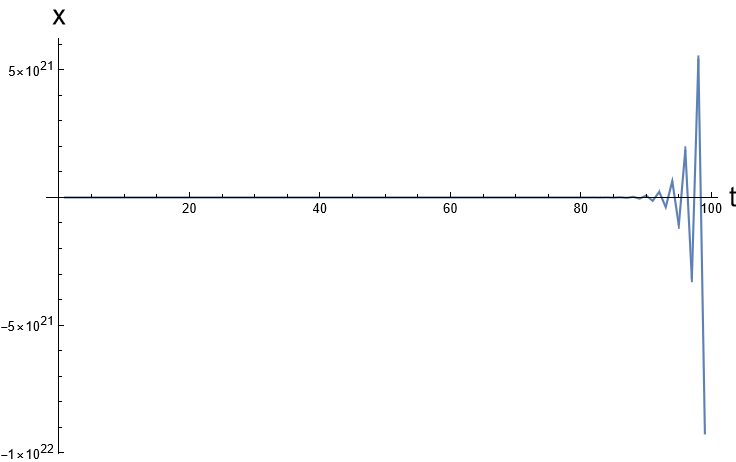}
        \label{jcase2d}
        } \hspace{0.3cm}
    \subfigure[$a=-1.0874, b=0.8667+0.04332i$]{
\includegraphics[height=1.5in,width=1.5in,keepaspectratio]{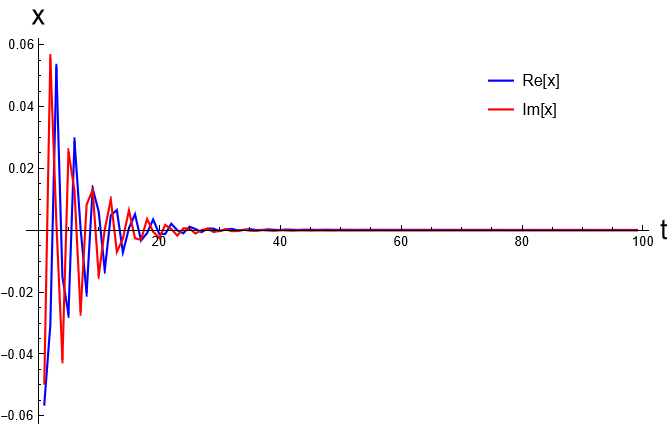}
        \label{kcase2a}
    } \hspace{0.3cm}
    \subfigure[$a=-1.0874, b=0.8424-0.04686i$]{
\includegraphics[height=1.5in,width=1.5in,keepaspectratio]{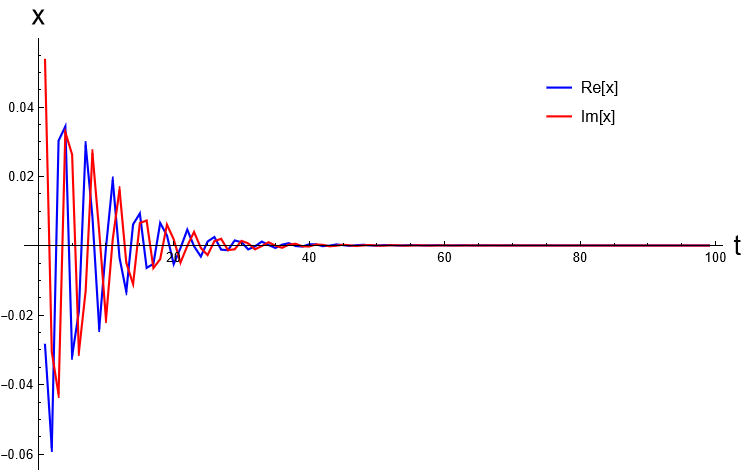}
        \label{kcase2b}
    } \hspace{0.3cm}
    \subfigure[$a=-1.0874, b=0.7+0.02i$]{
\includegraphics[height=1.5in,width=1.5in,keepaspectratio]{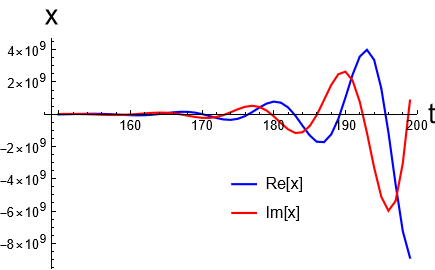}
        \label{kcase2c}
    } \hspace{0.3cm}
    \subfigure[$a=-1.0874, b=0.75-0.04i$]{
\includegraphics[height=1.5in,width=1.5in,keepaspectratio]{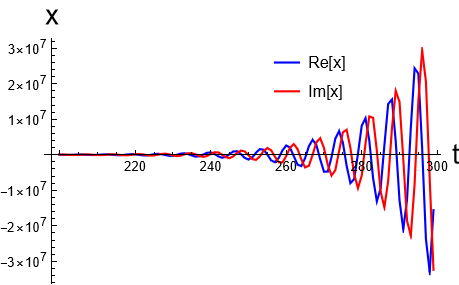}
        \label{kcase2d}
    } \hspace{0.3cm}
    \subfigure[$a=-1.0874, \, b=0.95$]{
\includegraphics[height=1.5in,width=1.5in,keepaspectratio]{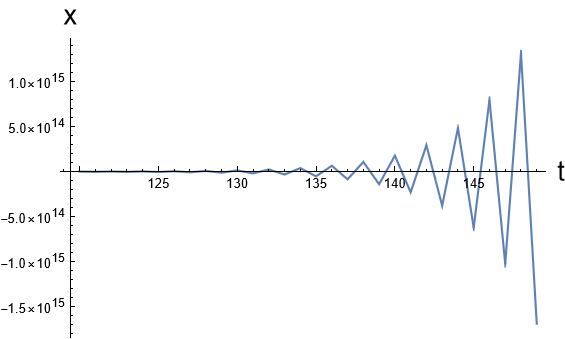}
        \label{kcase2e}
    } \hspace{0.3cm}
    \subfigure[$a=-1.1487, b=0.8325+0.07685i$]{
\includegraphics[height=1.5in,width=1.5in,keepaspectratio]{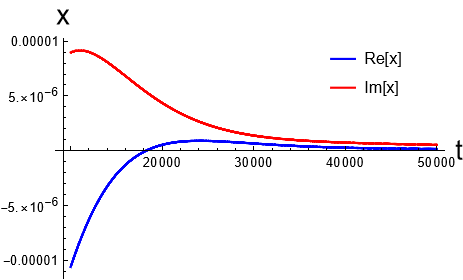}
        \label{mcase2a}
        } \hspace{0.3cm}
    \subfigure[$a=-1.1487, b=0.8616-0.0896i$]{
\includegraphics[height=1.5in,width=1.5in,keepaspectratio]{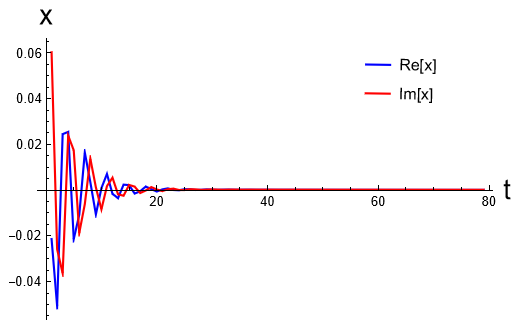}
        \label{mcase2b}
        } \hspace{0.3cm}
    \subfigure[$a=-1.1487, b=0.6458+0.03003i$]{
\includegraphics[height=1.5in,width=1.5in,keepaspectratio]{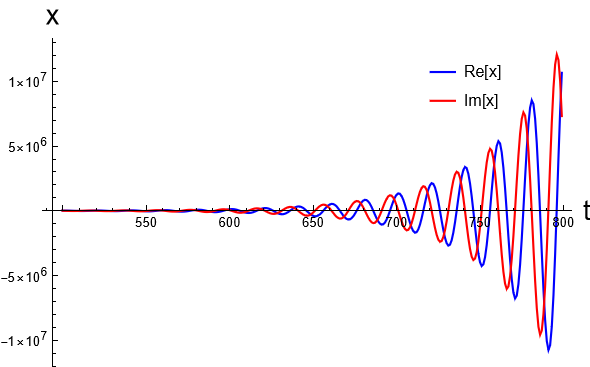}
        \label{mcase2c}
        } \hspace{0.3cm}
    \subfigure[$a=-1.1487, b=0.6497-0.03331i$]{
\includegraphics[height=1.5in,width=1.5in,keepaspectratio]{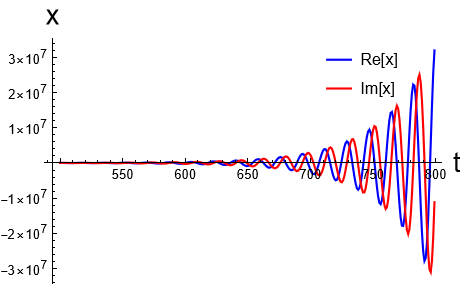}
        \label{mcase2d}
    }
    \caption{ The solutions of the system (\ref{1}) for $\alpha=0.4$, $\beta=0.2$ and different values of $a$ and $b$.}
    \label{1c1}
\end{figure}

%%%%%%%%%%%%%%%%%%%%%%%%%%%%%%%%%%%%%%%%%%%%%%%%%%%%%%%%%%%%%%%%%%%%%%%%%%%%%%%%%%%%%%%%%%%%%%%%%%
\section{Conclusion} \label{con.}
We studied the two-term linear fractional difference equation $
\Delta^{\alpha} x(t) + a \, \Delta^{\beta} x(t+\alpha-\beta) =(b-1)x(t+\alpha-1)$. We investigated the equivalent sequence form of the initial value problem (\ref{1}) with $0<\beta<\alpha \leq 1$ and validated it. We derived the expression for the boundary curve $\gamma(\theta)$ and also sketched the stable regions for the different sets of three parameters $a$, $\alpha$, and $\beta$ in a complex plane for the system. Additionally, we discussed the bifurcations for $a<0$ with the cases $0.5<\beta<\alpha<1$ and $0<\beta<\alpha<0.5$ in detail. We observed the rich dynamics in the considered cases. We provided numerical experiments to support our analytical results. 	
	
\section*{Acknowledgments}
 Ch. Janardhan thanks the University Grants Commission, New Delhi, India, for financial support (No. F.14-34/2011(CPP-II)).
	
\section*{Data Availability }
The dataset referenced in the context of this research work is stored at the link \url{https://github.com/Janardhan3233/2term}.

\bibliographystyle{plain}      % 
\bibliography{ref.bib}

\begin{thebibliography}{50}
\providecommand{\natexlab}[1]{#1}
\providecommand{\url}[1]{\texttt{#1}}
\expandafter\ifx\csname urlstyle\endcsname\relax
  \providecommand{\doi}[1]{doi: #1}\else
  \providecommand{\doi}{doi: \begingroup \urlstyle{rm}\Url}\fi

\bibitem[Aris(1994)]{aris1994mathematical}
Rutherford Aris.
\newblock \emph{Mathematical modelling techniques}.
\newblock Courier Corporation, 1994.

\bibitem[Berry and Houston(1995)]{berry1995mathematical}
John Berry and Ken Houston.
\newblock \emph{Mathematical modelling}.
\newblock Butterworth-Heinemann, 1995.

\bibitem[Gombert and Nielsen(2000)]{gombert2000mathematical}
Andreas~Karoly Gombert and Jens Nielsen.
\newblock Mathematical modelling of metabolism.
\newblock \emph{Current opinion in biotechnology}, 11\penalty0 (2):\penalty0
  180--186, 2000.

\bibitem[Elaydi(2005{\natexlab{a}})]{elaydi2005introduction}
Saber Elaydi.
\newblock \emph{An introduction to difference equations}.
\newblock Springer, 2005{\natexlab{a}}.

\bibitem[May(1976)]{may1976simple}
Robert~M May.
\newblock Simple mathematical models with very complicated dynamics.
\newblock \emph{Nature}, 261\penalty0 (5560):\penalty0 459--467, 1976.

\bibitem[Joshi et~al.(2024)Joshi, Bhalekar, and Gade]{joshi2024stability}
Divya~D Joshi, Sachin Bhalekar, and Prashant~M Gade.
\newblock Stability analysis of fractional difference equations with delay.
\newblock \emph{Chaos: An Interdisciplinary Journal of Nonlinear Science},
  34\penalty0 (5), 2024.

\bibitem[Kotulski et~al.(1999)Kotulski, Szczepa{\'n}ski, G{\'o}rski,
  Paszkiewicz, and Zugaj]{kotulski1999application}
Zbigniew Kotulski, Janusz Szczepa{\'n}ski, Karol G{\'o}rski, Andrzej
  Paszkiewicz, and Anna Zugaj.
\newblock Application of discrete chaotic dynamical systems in
  cryptography—dcc method.
\newblock \emph{International Journal of Bifurcation and Chaos}, 9\penalty0
  (06):\penalty0 1121--1135, 1999.

\bibitem[Singh and Sinha(2010)]{singh2010chaos}
Narendra Singh and Aloka Sinha.
\newblock Chaos-based secure communication system using logistic map.
\newblock \emph{Optics and Lasers in Engineering}, 48\penalty0 (3):\penalty0
  398--404, 2010.

\bibitem[Zhang and Cao(2014)]{zhang2014novel}
Xianhan Zhang and Yang Cao.
\newblock A novel chaotic map and an improved chaos-based image encryption
  scheme.
\newblock \emph{The Scientific World Journal}, 2014\penalty0 (1):\penalty0
  713541, 2014.

\bibitem[Khrennikov(2004)]{khrennikov2004p}
A~Khrennikov.
\newblock p-adic discrete dynamical systems and their applications in physics
  and cognitive sciences.
\newblock \emph{Russian Journal of Mathematical Physics}, 11\penalty0
  (1):\penalty0 45--70, 2004.

\bibitem[Huang and Wang(2015)]{huang2015applications}
Chun~Miao Huang and Wei~Ping Wang.
\newblock Applications of difference equation in population forecasting model.
\newblock \emph{Advanced Materials Research}, 1079:\penalty0 664--667, 2015.

\bibitem[Sandefur(1990)]{sandefur1990discrete}
James~T Sandefur.
\newblock \emph{Discrete dynamical systems: Theory and applications}.
\newblock Clarendon Press, 1990.

\bibitem[Tu(2012)]{tu2012dynamical}
Pierre~NV Tu.
\newblock \emph{Dynamical systems: an introduction with applications in
  economics and biology}.
\newblock Springer Science \& Business Media, 2012.

\bibitem[Bahi and Guyeux(2013)]{bahi2013discrete}
Jacques~M Bahi and Christophe Guyeux.
\newblock \emph{Discrete dynamical systems and chaotic machines: theory and
  applications}.
\newblock CRC Press, 2013.

\bibitem[AlSharawi et~al.(2014)AlSharawi, Cushing, and
  Elaydi]{alsharawi2014theory}
Ziyad AlSharawi, Jim~M Cushing, and Saber Elaydi.
\newblock Theory and applications of difference equations and discrete
  dynamical systems.
\newblock \emph{Springer Proceedings in Mathematics \& Statistics}, 102, 2014.

\bibitem[Podlubny(1998)]{podlubny1998fractional}
Igor Podlubny.
\newblock \emph{Fractional differential equations: an introduction to
  fractional derivatives, fractional differential equations, to methods of
  their solution and some of their applications}, volume 198.
\newblock elsevier, 1998.

\bibitem[Kilbas et~al.(2006)Kilbas, Srivastava, and Trujillo]{kilbas2006theory}
A~Anatolii~Aleksandrovich Kilbas, Hari~Mohan Srivastava, and Juan~J Trujillo.
\newblock \emph{Theory And Applications of Fractional Differential Equations},
  volume 204.
\newblock Elsevier Science Limited, 2006.

\bibitem[Diethelm and Ford(2002)]{diethelm2002analysis}
Kai Diethelm and Neville~J Ford.
\newblock Analysis of fractional differential equations.
\newblock \emph{Journal of Mathematical Analysis and Applications},
  265\penalty0 (2):\penalty0 229--248, 2002.

\bibitem[Daftardar-Gejji and Babakhani(2004)]{daftardar2004analysis}
Varsha Daftardar-Gejji and Azizollah Babakhani.
\newblock Analysis of a system of fractional differential equations.
\newblock \emph{Journal of Mathematical Analysis and Applications},
  293\penalty0 (2):\penalty0 511--522, 2004.

\bibitem[Bhalekar(2016)]{bhalekar2016stability}
Sachin Bhalekar.
\newblock Stability and bifurcation analysis of a generalized scalar delay
  differential equation.
\newblock \emph{Chaos: An Interdisciplinary Journal of Nonlinear Science},
  26\penalty0 (8), 2016.

\bibitem[Gupta and Bhalekar(2024)]{gupta2024fractional}
Deepa Gupta and Sachin Bhalekar.
\newblock Fractional order sunflower equation: Stability, bifurcation and
  chaos.
\newblock \emph{The European Physical Journal Special Topics}, pages 1--11,
  2024.

\bibitem[Bhalekar and Dutta(2025)]{bhalekar2025analysis}
Sachin Bhalekar and Pragati Dutta.
\newblock Analysis of a class of two-delay fractional differential equation.
\newblock \emph{Chaos: An Interdisciplinary Journal of Nonlinear Science},
  35\penalty0 (1), 2025.

\bibitem[Magin(2012)]{magin2012fractional}
Richard~L Magin.
\newblock Fractional calculus in bioengineering: A tool to model complex
  dynamics.
\newblock In \emph{Proceedings of the 13th International Carpathian Control
  Conference (ICCC)}, pages 464--469. IEEE, 2012.

\bibitem[Mainardi(2022)]{mainardi2022fractional}
Francesco Mainardi.
\newblock \emph{Fractional calculus and waves in linear viscoelasticity: an
  introduction to mathematical models}.
\newblock World Scientific, 2022.

\bibitem[Hilfer(2000)]{hilfer2000applications}
Rudolf Hilfer.
\newblock \emph{Applications of fractional calculus in physics}.
\newblock World scientific, 2000.

\bibitem[Lubich(1986)]{lubich1986discretized}
Ch~Lubich.
\newblock Discretized fractional calculus.
\newblock \emph{SIAM Journal on Mathematical Analysis}, 17\penalty0
  (3):\penalty0 704--719, 1986.

\bibitem[Miller and Ross(1989)]{miller1989fractional}
Kenneth~S Miller and Bertram Ross.
\newblock Fractional difference calculus.
\newblock In \emph{Proceedings of the international symposium on univalent
  functions, fractional calculus and their applications}, pages 139--152, 1989.

\bibitem[Atici and Eloe(2009)]{atici2009initial}
Ferhan Atici and Paul Eloe.
\newblock Initial value problems in discrete fractional calculus.
\newblock \emph{Proceedings of the American mathematical society}, 137\penalty0
  (3):\penalty0 981--989, 2009.

\bibitem[Fulai et~al.(2011)Fulai, Xiannan, and Yong]{fulai2011existence}
Chen Fulai, Luo Xiannan, and Zhou Yong.
\newblock Existence results for nonlinear fractional difference equation.
\newblock \emph{Advances in Difference Equations}, 12, 2011.

\bibitem[Goodrich and Peterson(2015)]{goodrich2015discrete}
Christopher Goodrich and Allan~C Peterson.
\newblock \emph{Discrete fractional calculus}, volume~10.
\newblock Springer, 2015.

\bibitem[Ferreira et~al.(2022)]{ferreira2022discrete}
Rui~AC Ferreira et~al.
\newblock \emph{Discrete fractional calculus and fractional difference
  equations}.
\newblock Springer, 2022.

\bibitem[Abu-Saris and Al-Mdallal(2013)]{abu2013asymptotic}
Raghib Abu-Saris and Qasem Al-Mdallal.
\newblock On the asymptotic stability of linear system of fractional-order
  difference equations.
\newblock \emph{Fractional Calculus and Applied Analysis}, 16\penalty0
  (3):\penalty0 613--629, 2013.

\bibitem[{\v{C}}erm{\'a}k et~al.(2015){\v{C}}erm{\'a}k, Gy{\H{o}}ri, and
  Nechv{\'a}tal]{vcermak2015explicit}
Jan {\v{C}}erm{\'a}k, Istv{\'a}n Gy{\H{o}}ri, and Lud{\u{e}}k Nechv{\'a}tal.
\newblock On explicit stability conditions for a linear fractional difference
  system.
\newblock \emph{Fractional Calculus and Applied Analysis}, 18\penalty0
  (3):\penalty0 651--672, 2015.

\bibitem[Bhalekar et~al.(2022)Bhalekar, Gade, and
  Joshi]{bhalekar2022stabilitycomplexorder}
Sachin Bhalekar, Prashant~M Gade, and Divya Joshi.
\newblock Stability and dynamics of complex order fractional difference
  equations.
\newblock \emph{Chaos, Solitons \& Fractals}, 158:\penalty0 112063, 2022.

\bibitem[Bhalekar and Gade(2022)]{bhalekar2022stabilitycoupledmap}
Sachin Bhalekar and Prashant~M Gade.
\newblock Stability analysis of fixed point of fractional-order coupled map
  lattices.
\newblock \emph{Communications in Nonlinear Science and Numerical Simulation},
  113:\penalty0 106587, 2022.

\bibitem[Bhalekar and Gade(2023)]{bhalekar2023fractional}
Sachin Bhalekar and Prashant~M Gade.
\newblock Fractional-order periodic maps: Stability analysis and application to
  the periodic-2 limit cycles in the nonlinear systems.
\newblock \emph{Journal of Nonlinear Science}, 33\penalty0 (6):\penalty0 119,
  2023.

\bibitem[Joshi et~al.(2023)Joshi, Bhalekar, and Gade]{joshi2023controlling}
Divya~D Joshi, Sachin Bhalekar, and Prashant~M Gade.
\newblock Controlling fractional difference equations using feedback.
\newblock \emph{Chaos, Solitons \& Fractals}, 170:\penalty0 113401, 2023.

\bibitem[Bhalekar et~al.(2025)Bhalekar, Chevala, and
  Gade]{bhalekar2025dynamical}
Sachin Bhalekar, Janardhan Chevala, and Prashant~M Gade.
\newblock Dynamical analysis of fractional order generalized logistic map.
\newblock \emph{Computational Mathematics and Mathematical Physics},
  65\penalty0 (2):\penalty0 424--441, 2025.

\bibitem[Daftardar-Gejji and
  Bhalekar(2008{\natexlab{a}})]{daftardar2008boundary}
Varsha Daftardar-Gejji and Sachin Bhalekar.
\newblock Boundary value problems for multi-term fractional differential
  equations.
\newblock \emph{Journal of Mathematical Analysis and Applications},
  345\penalty0 (2):\penalty0 754--765, 2008{\natexlab{a}}.

\bibitem[Daftardar-Gejji and
  Bhalekar(2008{\natexlab{b}})]{daftardar2008solving}
Varsha Daftardar-Gejji and Sachin Bhalekar.
\newblock Solving multi-term linear and non-linear diffusion--wave equations of
  fractional order by adomian decomposition method.
\newblock \emph{Applied Mathematics and Computation}, 202\penalty0
  (1):\penalty0 113--120, 2008{\natexlab{b}}.

\bibitem[{\v{C}}erm{\'a}k and Kisela(2015{\natexlab{a}})]{vcermak2015stability}
Jan {\v{C}}erm{\'a}k and Tom{\'a}{\v{s}} Kisela.
\newblock Stability properties of two-term fractional differential equations.
\newblock \emph{Nonlinear Dynamics}, 80\penalty0 (4):\penalty0 1673--1684,
  2015{\natexlab{a}}.

\bibitem[Bhalekar and Gupta(2024)]{bhalekar2024stability}
Sachin Bhalekar and Deepa Gupta.
\newblock Stability and bifurcation analysis of two-term fractional
  differential equation with delay.
\newblock \emph{arXiv preprint arXiv:2404.01824}, 2024.

\bibitem[{\v{C}}erm{\'a}k and
  Kisela(2015{\natexlab{b}})]{vcermak2015asymptotic}
Jan {\v{C}}erm{\'a}k and Tom{\'a}{\v{s}} Kisela.
\newblock Asymptotic stability of dynamic equations with two fractional terms:
  continuous versus discrete case.
\newblock \emph{Fractional Calculus and Applied Analysis}, 18\penalty0
  (2):\penalty0 437--458, 2015{\natexlab{b}}.

\bibitem[Srivastava and Rai(2010)]{srivastava2010multi}
Vineet Srivastava and KN~Rai.
\newblock A multi-term fractional diffusion equation for oxygen delivery
  through a capillary to tissues.
\newblock \emph{Mathematical and Computer Modelling}, 51\penalty0
  (5-6):\penalty0 616--624, 2010.

\bibitem[Bastos et~al.(2011)Bastos, Ferreira, and Torres]{bastos2011discrete}
Nuno~RO Bastos, Rui~AC Ferreira, and Delfim~FM Torres.
\newblock Discrete-time fractional variational problems.
\newblock \emph{Signal Processing}, 91\penalty0 (3):\penalty0 513--524, 2011.

\bibitem[Ferreira and Torres(2011)]{ferreira2011fractional}
Rui~AC Ferreira and Delfim~FM Torres.
\newblock Fractional h-difference equations arising from the calculus of
  variations.
\newblock \emph{Applicable Analysis and Discrete Mathematics}, pages 110--121,
  2011.

\bibitem[Mozyrska et~al.(2015)Mozyrska, Wyrwas, et~al.]{mozyrska2015transform}
Dorota Mozyrska, Ma{\l}gorzata Wyrwas, et~al.
\newblock The-transform method and delta type fractional difference operators.
\newblock \emph{Discrete Dynamics in Nature and Society}, 2015, 2015.

\bibitem[Atici and Eloe(2007)]{atici2007transform}
Ferhan~M Atici and Paul~W Eloe.
\newblock A transform method in discrete fractional calculus.
\newblock \emph{International Journal of Difference Equations}, 2\penalty0 (2),
  2007.

\bibitem[Elaydi(2005{\natexlab{b}})]{elaydi2005systems}
Saber Elaydi.
\newblock Systems of linear difference equations.
\newblock \emph{An Introduction to Difference equations}, pages 117--172,
  2005{\natexlab{b}}.

\bibitem[Hirsch et~al.(2013)Hirsch, Smale, and Devaney]{hirsch2013differential}
Morris~W Hirsch, Stephen Smale, and Robert~L Devaney.
\newblock \emph{Differential equations, dynamical systems, and an introduction
  to chaos}.
\newblock Academic press, 2013.

\end{thebibliography}

\end{document}